\documentclass[11pt]{article}

\usepackage{booktabs}
\usepackage[margin=15pt,font=small,labelfont={bf,sf},justification=justified]{caption}
\usepackage[superscript,biblabel]{cite}
\usepackage{color}
\usepackage{float}
\usepackage{graphicx}
\usepackage{lscape}
\usepackage[bbgreekl]{mathbbol}
\usepackage{mathtools}
\usepackage{multibib}
\usepackage{multirow}
\usepackage{tabularx}
\usepackage{titleref}
\usepackage{url}

\usepackage{tikz}
\usepackage{amsmath}

\usetikzlibrary{fadings}
\usetikzlibrary{patterns}
\usetikzlibrary{shadows.blur}
\usetikzlibrary{shapes}

\newif\ifbembo
\newif\ifcharter
\newif\iferewhon
\newif\iflibertine
\newif\iflibertinealt
\newif\ifpalantino
\newif\iftimesnewroman

\bembofalse
\chartertrue
\erewhonfalse
\libertinefalse
\libertinealtfalse
\palantinofalse
\timesnewromanfalse

\ifbembo
  \usepackage[p,osf]{ETbb} 
  \usepackage[scaled=.95,type1]{cabin} 
  \usepackage[varqu,varl]{zi4}
  \usepackage[T1]{fontenc}
  \usepackage[libertine,vvarbb]{newtxmath}
  \usepackage[cal=boondoxo,bb=boondox]{mathalfa}
\fi

\ifcharter
  \usepackage[scaled=.98,p]{XCharter}
  \usepackage[scaled=1.04,varqu,varl]{inconsolata}
  \usepackage[type1]{cabin}
  \usepackage[uprightscript,xcharter,vvarbb,scaled=1.05]{newtxmath}
  \usepackage[cal=euler,scr=boondoxo,bb=boondox]{mathalpha}  
\fi

\iferewhon
  \usepackage[p,osf,scaled=.98,space]{erewhon} 
  \usepackage[varqu,varl]{inconsolata} 
  \usepackage[type1,scaled=.95]{cabin} 
  \usepackage[utopia,vvarbb]{newtxmath}
\fi

\iflibertine
  \usepackage[sb]{libertinus}
  \usepackage[T1]{fontenc}
  \usepackage{textcomp}
  \usepackage[varqu,varl]{zi4}
  \usepackage{libertinust1math} 
  \usepackage[scr=boondoxo,bb=boondox]{mathalpha} 
\fi

\iflibertinealt
  \usepackage[sb]{libertinus}
  \usepackage[T1]{fontenc}
  \usepackage{textcomp}
  \usepackage{libertinust1math} 
  \usepackage[scr=boondoxo,bb=boondox]{mathalpha} 
\fi

\ifpalantino
  \usepackage[largesc]{newpxtext}
  \usepackage[vvarbb]{newpxmath}
\fi

\iftimesnewroman
  \usepackage{newtxtext} 
  \usepackage[scaled=0.95]{inconsolata} 
  \usepackage[vvarbb]{newtxmath} 
\fi

\usepackage[T1]{fontenc}
\usepackage[final,protrusion=true,expansion=true]{microtype}

\usepackage{titling}
\usepackage{authblk} 
\pretitle{\begin{center}\bfseries\sffamily\LARGE}
\posttitle{\end{center}}

\predate{\begin{center}\normalsize}
\postdate{\end{center}}

\usepackage{abstract}

\usepackage{sectsty} 
\allsectionsfont{\sffamily}

\usepackage[left=1in,right=1in,top=1in,bottom=1in,includefoot,heightrounded]{geometry}

\newcites{supp}{SUPPLEMENTAL REFERENCES}

\makeatletter
\patchcmd{\LS@rot}{90}{-90}{}{}
\patchcmd{\endlandscape}{90}{-90}{}{}
\makeatother

\newcommand*{\tran}{^{\mkern-1.5mu\mathsf{T}}}
\newcommand{\cc}{{\text{c}}}
\newcommand{\Lmax}{{\text{L}_\text{max}}}
\newcommand{\Lrms}{{\text{L}_\text{rms}}}
\newcommand{\Davg}{{\text{D}_\text{avg}}}
\newcommand{\Drms}{{\text{D}_\text{rms}}}

\begin{document}

\title{An Immersed Interface Method for Incompressible Flows and Geometries with Sharp Features}

\author[1]{Michael J. Facci}
\author[2]{Ebrahim M.~Kolahdouz}
\author[1,3--6]{Boyce E.~Griffith}
\affil[1]{Department of Mathematics, University of North Carolina, Chapel Hill, NC, USA}
\affil[2]{Flatiron Institute, Simons Foundation, New York, NY, USA}
\affil[3]{Department of Biomedical Engineering, University of North Carolina, Chapel Hill, NC, USA}
\affil[4]{Carolina Center for Interdisciplinary Applied Mathematics, University of North Carolina, Chapel Hill, NC, USA}
\affil[5]{Computational Medicine Program, University of North Carolina School of Medicine, Chapel Hill, NC, USA}
\affil[6]{McAllister Heart Institute, University of North Carolina School of Medicine, Chapel Hill, NC, USA\vspace{\baselineskip}}

\maketitle
\begin{abstract}
   The immersed interface method (IIM) for models of fluid flow and fluid-structure interaction imposes jump conditions that capture stress discontinuities generated by forces that are concentrated along immersed boundaries. Most prior work using the IIM for fluid dynamics applications has focused on smooth interfaces, but boundaries with sharp features such as corners and edges can appear in practical analyses, particularly on engineered structures. The present study builds on our work to integrate finite element-type representations of interface geometries with the IIM. Initial realizations of this approach used a continuous Galerkin (CG) finite element discretization for the boundary, but as we show herein, these approaches generate large errors near sharp geometrical features. To overcome this difficulty, this study introduces an IIM approach using a discontinuous Galerkin (DG) representation of the jump conditions. Numerical examples explore the impacts of different interface representations on accuracy for both smooth and sharp boundaries, particularly flows interacting with fixed interface configurations. We demonstrate that using a DG approach provides accuracy that is comparable to the CG method for smooth cases. Further, we identify a time step size restriction for the CG representation that is directly related to the sharpness of the geometry. In contrast, time step size restrictions imposed by DG representations are demonstrated to be nearly insensitive to the presence of sharp features.
\end{abstract}


\section{Introduction}

Simulating the dynamics of interactions between fluid flows and immersed boundaries has a broad range of applications in science and engineering. These include heart valves \cite{griffith_immersed_2012,peskin_flow_1972}, bio-locomotion for swimmers \cite{zeng_immersed_2023}, and inferior vena cava filters \cite{kolahdouz_sharp_2023}. The immersed interface method (IIM) for incompressible flows imposes stress discontinuities that arise where interfacial forces are concentrated along a boundary that is immersed in a fluid.
The IIM originated in work by Mayo{\cite{mayo_fast_1984}} and by LeVeque and Li{\cite{leveque_immersed_1994}} for the solution of
elliptic partial differential equations with singular forces. The IIM was
extended to the two-dimensional solution of the incompressible Navier-Stokes
equations by Lee and LeVeque\cite{lee_immersed_2003} and Lai and Li \cite{li_immersed_2001}. 
Jump conditions for velocity, pressure, and their derivatives for fluid flows in three spatial dimensions were derived by Xu and Wang \cite{xu_systematic_2006}.
Immersed interface methods developed by Thekkethil and Sharma \cite{thekkethil_level_2019} and Xu et al. \cite{xu_level-set_2020} use level set representations of the
geometry. In contrast,
the IIM described by Kolahdouz {\cite{kolahdouz_immersed_2020,kolahdouz_sharp_2021,kolahdouz_sharp_2023}} deals with
interfaces represented as triangulated surfaces. Such surface representations can be readily use with finite element structural models. In principle, this allows
for geometries with sharp features to be represented, such as a bileaflet mechanical
heart valve {\cite{kolahdouz_sharp_2021}}. However, most prior work using the IIM for fluid dynamics and FSI has focused on intrinsically smooth interfaces, in which the limit surface is $C^\infty$ even if the discrete representation is only $C^0$. One exception is the work by Le et. al \cite{LE2006109}, which presents an IIM handling both smooth and sharp geometry with an arc length parameterization of the interface. 
Another work by Liu and Xu\cite{sharp_liu_xu} presents an IIM for sharp geometries using polygonal curves. A limitation of both studies is their parameterization of the interface, which poorly suited for complex engineering geometries, particularly in three spatial dimensions.

This study considers the performance of our immersed interface method for discrete geometries in the context of interfaces with sharp features such as corners or edges. 
These types of features occur in engineering applications, such as mechanical heart valves, turbine blades, and airfoils. 
We revisit our original approach of the discrete approximation of the geometry and jump conditions. The lowest order jump conditions are computed using the normal and tangential components of the interfacial forces. 
If the interface is $C^0$ but not $C^1$, then the jump conditions themselves are also discontinuous at element junctions. 
Jump conditions are defined in the continuous equations in terms of the surface normal, and are computed in our immersed interface methodology using the interfacial surface element's normal vector, which is
discontinuous between elements if the geometry is captured via $C^0$ Lagrange basis functions. In the original work, we used an $L^2$ projection of the space
spanned by the standard continuous Lagrangian basis functions to construct continuous
jump condition values along the $C^0$ representation of smooth interfaces and
used Gaussian quadrature to approximate the numerical integration. However, projecting the jump conditions onto
a continuous basis produces an $O (1)$ error at sharp corners. One approach to ameliorate these errors is to regularize the sharp corners or edges by
introducing a fillet with a small radius for the fillet arc. Closely following
the toplogical changes around sharp corners with a piecewise parametric
surface mesh, however, inevitably requires introducing elements around the
corners and edges with sizes that are orders of magnitudes smaller than the
average element size in the bulk mesh. Resolving these sliver elements is challenging, particularly in a fluid-structure interaction (FSI)
framework for which Eulerian grid resolution needs to be comparable to the
Lagrangian element size.

This paper introduces a discontinuous Galerkin{\cite{cockburn_discontinuous_2017}} immersed interface method (DG-IIM) that uses an intrinsically discontinuous representation of the jump conditions to improve the accuracy of the IIM for interfaces with sharp features, without smoothing edges or additional mesh refinement near
these singular edges and corners. Components of the force are projected onto discontinuous Lagrange functions with element-local support. The discontinuous Lagrange basis functions are
similar to the Lagrange family but without interelement continuity. Each node
contains extra degrees of freedom for each element it shares. We show that for sharp geometries, this formulation more accurately
imposes prescribed stationary rigid body motion. Numerical experiments reveal that
discontinuous elements also allow for a larger time step size, increasing efficiency
of numerical simulations.
We treat the interfaces with this representation during the force spreading step, in which discontinuous projections of the jump conditions are used to compute correction terms in the Eulerian finite difference stencils. The key contribution of this study is that it demonstrates that this DG-IIM proves to be an effective coupling scheme for solving FSI with
sharp edges and vertices. It accurately reproduces literature results for
smooth geometry, while also achieving significantly more accuracy and efficiency than a continuous Galerkin immersed interface method (CG-IIM) formulation, which projects the jump conditions onto standard $C^0$ Lagrangian basis functions.
Numerical experiments show that increasing acuteness in
interfacial edges does not affect the maximum stable time step size for the DG-IIM, but decreases substantially for a CG-IIM
formulation. 

\section{Continuous Equations of Motion}

This section outlines the equations of motion for a viscous incompressible fluid with interfacial forces concentrated on an immersed boundary. This boundary both applies force to the fluid and moves with the local fluid velocity. The fluid domain, $\Omega $, is partitioned
into two subregions, $\Omega_t^+$ (exterior) and $\Omega_t^-$ (interior), by an
interface $\Gamma_t = \overline{\Omega_t^+} \bigcap \overline{\Omega_t^-}$.  We use the initial coordinates $\mathbf{X}$ of the boundary as Lagrangian coordinates, so that $\mathbf{X} \in \Gamma_0$, with corresponding current coordinates $\boldsymbol{\chi}(\mathbf{X},t) \in \Gamma_t$. 
The equations of motion are
\begin{align}
  \rho \frac{\mathrm{D} \mathbf{u}(\mathbf{x},t)}{\mathrm{D}t}& = - \nabla p (\mathbf{x}, t) + \mu
  \nabla^2 \mathbf{u} (\mathbf{x}, t), & \mathbf{x} \in \Omega,  \\
  \nabla \cdot \mathbf{u} (\mathbf{x}, t) & = 0, & \mathbf{x} \in
  \Omega, \\
  \llbracket \mathbf{u} (\mathbf{\boldsymbol{\chi}} (\mathbf{X}, t), t) \rrbracket &
  = 0,  &\mathbf{X} \in \Gamma_0,  \\
  \llbracket p (\mathbf{\boldsymbol{\chi}} (\mathbf{X}, t), t) \rrbracket & = -\jmath^{-
  1} (\mathbf{X}, t)\, \mathbf{F} (\mathbf{X}, t) \cdot \mathbf{n} (\mathbf{X}, t), & \mathbf{X} \in \Gamma_0, \\
  \mu\left\llbracket  \frac{\partial \mathbf{u} (\mathbf{\boldsymbol{\chi}}
  (\mathbf{X}, t), t)}{\partial \mathbf{x}_i} \right\rrbracket &
  =   (\mathbb{I}- \mathbf{n} (\mathbf{X}, t)
  \mathbf{n} (\mathbf{X}, t)\tran) \jmath^{- 1} (\mathbf{X}, t)\, \mathbf{F} (\mathbf{X}, t) n_i  (\mathbf{X}, t), &
  \mathbf{X} \in \Gamma_0, 
  \\  \mathbf{U} (\mathbf{X}, t)  &
  = \mathbf{u} (\boldsymbol{\chi}(\mathbf{X}, t), t),  &\mathbf{X} \in \Gamma_0, 
  \\ \frac{\partial \boldsymbol{\chi(\mathbf{X},t)}}{\partial t}     &
  = \mathbf{U}(\mathbf{X},t),  &\mathbf{X} \in \Gamma_0, 
\end{align}
in which $\rho$, $\mu$, $\mathbf{u} (\mathbf{x}, t)$, and $p
(\mathbf{x}, t)$ are the fluid's mass
density, dynamic viscosity, velocity, and pressure. 
$\mathbf{F} (\mathbf{X}, t)$ and $\mathbf{U} (\mathbf{X}, t)$ are the interfacial force and velocity, respectively, in Lagrangian material coordinates, and $\mathbf{n} (\mathbf{X}, t)$ is the surface normal vector in the interface's current configuration.
 $\jmath^{- 1} =
\frac{{\mathrm{d}a}}{{\mathrm{d}A}}$ relates surface area along the interface in the current $(\mathrm{d}a)$ and reference $(\mathrm{d}A)$ configurations. For stationary bodies, $\jmath^{- 1} = 1$.
The jump conditions express the discontinuity in fluid traction generated by the singular force along the interface. 
For a general function $\phi(\mathbf{x},t)$ and position $\mathbf{x}= \mathbf{\boldsymbol{\chi}}
  (\mathbf{X}, t)$ along the interface, the jump condition is
\begin{equation}
  \llbracket \phi (\mathbf{x}, t) \rrbracket \equiv \phi^+ (\mathbf{x}, t) - \phi^- (\mathbf{x}, t) \equiv \lim_{\varepsilon
  \rightarrow 0^+} \phi (\mathbf{x} + \varepsilon \mathbf{n}
  (\mathbf{X}, t), t) - \lim_{\varepsilon \rightarrow 0^-} \phi
  (\mathbf{x} - \varepsilon \mathbf{n} (\mathbf{X}, t), t).
\end{equation}
Equations (3)--(5) were derived in previous work by Lai and Li\cite{li_immersed_2001}, Xu and Wang\cite{xu_systematic_2006}, and Peskin and Printz\cite{PESKIN199333}. Both cases involve forces that are computed using a penalty method. In either case, we view the boundary as having a \textit{prescribed} configuration $\boldsymbol{\xi} (\mathbf{X},t)$ and an \textit{actual} configuration $\mathbf{\boldsymbol{\chi}} (\mathbf{X}, t)$.
Let $\mathbf{V}(\mathbf{X}, t) = \frac{\partial\boldsymbol{\xi} (\mathbf{X}, t)}{\partial t}$ be the velocity of the prescribed configuration. These two configurations are connected through penalty forces that are of the form

\begin{equation}
    \mathbf{F} (\mathbf{X}, t)  =  \kappa (\boldsymbol{\xi} (\mathbf{X},
    t) - \mathbf{\boldsymbol{\chi}} (\mathbf{X}, t)) + \eta \left( \mathbf{V}(\mathbf{X},
    t) - \mathbf{U}
    (\mathbf{X}, t) \right).
\end{equation}
$\kappa$ and $\eta$ are the penalty spring stiffness and damping coefficients. As $\kappa \rightarrow \infty$, this constraint exactly
imposes the motion\cite{IMFSI-griffith}. This study considers only the simplest case, in which the interface is stationary, so that $\boldsymbol{\chi}(\mathbf{X},t) = \mathbf{X}$ and $\mathbf{V}(\mathbf{X},t) = \bf0$.
The simplified force model becomes $\mathbf{F} (\mathbf{X}, t) = \kappa (\mathbf{X} - \mathbf{\boldsymbol{\chi}} (\mathbf{X}, t)) - \eta  \mathbf{U}
    (\mathbf{X}, t) $.

\section{Discrete Equations of Motion}

This section describes the spatial and temporal discretizations of the equations of motion. It also details the projection of the discrete jump conditions into both continuous and discontinous bases. We also define the discrete force spreading and velocity interpolation operators.
\subsection{Interface Representation}

We use a finite element representation of the immersed interface with
triangulation $\mathcal{T}_h$ of $\Gamma_0$, the reference configuration.
Consider elements $U_i$ such that $\mathcal{T}_h = \bigcup_i U_i$, with $i$ indexing
the mesh elements. The nodes of the mesh elements are $\{
\mathbf{X}_j \}_{j = 1}^M$ and have corresponding nodal (Lagrangian) basis functions $\{ \psi_j (\mathbf{X})
\}_{j = 1}^M$. Herein, we consider piecewise linear interface representations, in which $\psi_j$ is a piecewise linear Lagrange polynomial. The nodal basis functions $\psi_j(\mathbf{X})$ are continuous across elements but are not differentiable at the nodes. The current location of the interfacial
nodes at time $t$ are $\{ \mathbf{\boldsymbol{\chi}}_j (t) \}_{j = 1}^M$. In the finite
element space dictated by the subspace $S_h = \text{span} \{ \psi_j
(\mathbf{X}) \}_{j = 1}^M$, the configuration of the interface is
$ \mathbf{\boldsymbol{\chi}}_h (\mathbf{X}, t) = \sum_{j = 1}^M \mathbf{\boldsymbol{\chi}}_j (t)
  \psi_j (\mathbf{X})$.

\begin{figure}[t]
\center{}
\tikzset{every picture/.style={line width=0.75pt}} 

\begin{tikzpicture}[x=0.75pt,y=0.75pt,yscale=-1,xscale=1,scale = 1.2]

\draw  [line width=1.5]  (100,19.6) -- (301.2,19.6) -- (301.2,219.6) -- (100,219.6) -- cycle ;
\draw [line width=0.75]    (99.5,120.1) -- (301.2,120.1) ;
\draw [line width=0.75]    (100.5,170.1) -- (300.2,170.1) ;
\draw [line width=0.75]    (99.5,70.1) -- (301.2,70.1) ;
\draw [line width=0.75]    (150.5,18.6) -- (150.5,220.6) ;
\draw [line width=0.75]    (200.5,18.6) -- (200.5,220.6) ;
\draw [line width=0.75]    (249.5,19.6) -- (249.5,219.6) ;
\draw [color={rgb, 255:red, 74; green, 144; blue, 226 }  ,draw opacity=1 ]   (142.2,136.6) -- (159.2,219.6) ;
\draw [color={rgb, 255:red, 208; green, 2; blue, 27 }  ,draw opacity=1 ]   (178.2,76.6) -- (142.2,136.6) ;
\draw [color={rgb, 255:red, 245; green, 166; blue, 35 }  ,draw opacity=1 ]   (233.2,42.6) -- (178.2,76.6) ;
\draw [color={rgb, 255:red, 144; green, 19; blue, 254 }  ,draw opacity=1 ]   (302.2,35.6) -- (233.2,42.6) ;
\draw [color={rgb, 255:red, 74; green, 144; blue, 226 }  ,draw opacity=1 ]   (142.2,136.6) -- (116.12,144.05) ;
\draw [shift={(114.2,144.6)}, rotate = 344.05] [color={rgb, 255:red, 74; green, 144; blue, 226 }  ,draw opacity=1 ][line width=0.75]    (10.93,-3.29) .. controls (6.95,-1.4) and (3.31,-0.3) .. (0,0) .. controls (3.31,0.3) and (6.95,1.4) .. (10.93,3.29)   ;
\draw [color={rgb, 255:red, 208; green, 2; blue, 27 }  ,draw opacity=1 ]   (178.2,76.6) -- (155.9,62.66) ;
\draw [shift={(154.2,61.6)}, rotate = 32.01] [color={rgb, 255:red, 208; green, 2; blue, 27 }  ,draw opacity=1 ][line width=0.75]    (10.93,-3.29) .. controls (6.95,-1.4) and (3.31,-0.3) .. (0,0) .. controls (3.31,0.3) and (6.95,1.4) .. (10.93,3.29)   ;
\draw [color={rgb, 255:red, 245; green, 166; blue, 35 }  ,draw opacity=1 ]   (178.2,76.6) -- (164.29,55.28) ;
\draw [shift={(163.2,53.6)}, rotate = 56.89] [color={rgb, 255:red, 245; green, 166; blue, 35 }  ,draw opacity=1 ][line width=0.75]    (10.93,-3.29) .. controls (6.95,-1.4) and (3.31,-0.3) .. (0,0) .. controls (3.31,0.3) and (6.95,1.4) .. (10.93,3.29)   ;
\draw [color={rgb, 255:red, 245; green, 166; blue, 35 }  ,draw opacity=1 ]   (233.2,42.6) -- (222.16,22.36) ;
\draw [shift={(221.2,20.6)}, rotate = 61.39] [color={rgb, 255:red, 245; green, 166; blue, 35 }  ,draw opacity=1 ][line width=0.75]    (10.93,-3.29) .. controls (6.95,-1.4) and (3.31,-0.3) .. (0,0) .. controls (3.31,0.3) and (6.95,1.4) .. (10.93,3.29)   ;
\draw [color={rgb, 255:red, 208; green, 2; blue, 27 }  ,draw opacity=1 ]   (142.2,136.6) -- (119.9,122.66) ;
\draw [shift={(118.2,121.6)}, rotate = 32.01] [color={rgb, 255:red, 208; green, 2; blue, 27 }  ,draw opacity=1 ][line width=0.75]    (10.93,-3.29) .. controls (6.95,-1.4) and (3.31,-0.3) .. (0,0) .. controls (3.31,0.3) and (6.95,1.4) .. (10.93,3.29)   ;
\draw [color={rgb, 255:red, 144; green, 19; blue, 254 }  ,draw opacity=1 ]   (233.2,42.6) -- (229.53,20.57) ;
\draw [shift={(229.2,18.6)}, rotate = 80.54] [color={rgb, 255:red, 144; green, 19; blue, 254 }  ,draw opacity=1 ][line width=0.75]    (10.93,-3.29) .. controls (6.95,-1.4) and (3.31,-0.3) .. (0,0) .. controls (3.31,0.3) and (6.95,1.4) .. (10.93,3.29)   ;
\draw  [fill={rgb, 255:red, 0; green, 0; blue, 0 }  ,fill opacity=1 ] (140.2,136.6) .. controls (140.2,135.5) and (141.1,134.6) .. (142.2,134.6) .. controls (143.3,134.6) and (144.2,135.5) .. (144.2,136.6) .. controls (144.2,137.7) and (143.3,138.6) .. (142.2,138.6) .. controls (141.1,138.6) and (140.2,137.7) .. (140.2,136.6) -- cycle ;
\draw  [fill={rgb, 255:red, 0; green, 0; blue, 0 }  ,fill opacity=1 ] (176.34,76.6) .. controls (176.34,75.57) and (177.17,74.74) .. (178.2,74.74) .. controls (179.23,74.74) and (180.06,75.57) .. (180.06,76.6) .. controls (180.06,77.63) and (179.23,78.46) .. (178.2,78.46) .. controls (177.17,78.46) and (176.34,77.63) .. (176.34,76.6) -- cycle ;
\draw  [fill={rgb, 255:red, 0; green, 0; blue, 0 }  ,fill opacity=1 ] (231.34,42.6) .. controls (231.34,41.57) and (232.17,40.74) .. (233.2,40.74) .. controls (234.23,40.74) and (235.06,41.57) .. (235.06,42.6) .. controls (235.06,43.63) and (234.23,44.46) .. (233.2,44.46) .. controls (232.17,44.46) and (231.34,43.63) .. (231.34,42.6) -- cycle ;

\draw (114,34.4) node [anchor=north west][inner sep=0.75pt]    {$\Omega _{t}^{+}$};
\draw (213,138.4) node [anchor=north west][inner sep=0.75pt]    {$\Omega _{t}^{-}$};
\draw (106.5,145.5) node [anchor=north west][inner sep=0.75pt]  [font=\footnotesize]  {$\mathbf{n}(\mathbf{x} ,t)$};
\draw (153.5,41.5) node [anchor=north west][inner sep=0.75pt]  [font=\footnotesize]  {$\mathbf{n}(\mathbf{x} ,t)$};
\draw (214.5,50.5) node [anchor=north west][inner sep=0.75pt]  [font=\footnotesize]  {$\mathbf{n}(\mathbf{x} ,t)$};

\end{tikzpicture}
  \caption{Normal vectors are evaluated along the elements, but when directly obtained from the elements they are discontinuous at element junctions.}
\end{figure}
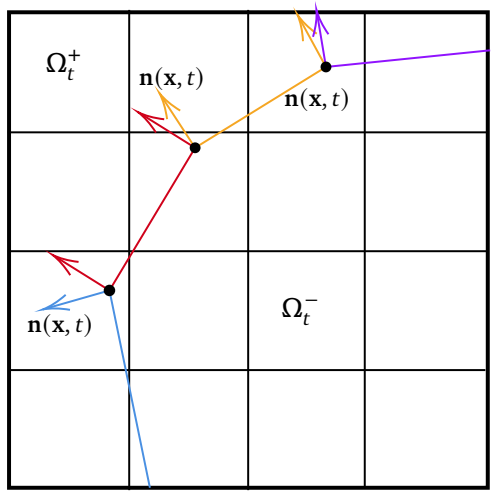

Directly evaluating the surface normal from a $C^0$ representation of the interface generates surface normals that are discontinuous at the junctions between elements (Figure 1). For a smooth surface, these discontinuities progressively vanish under grid refinement. In contrast, for surfaces with sharp features, these discontinuities persist. 

\subsection{Continuous and Discontinuous Projected Jump Conditions}
Jump conditions must be evaluated at arbitrary positions along the interface to evaluate the forcing terms associated with stencils that cut the interface. If using the normal vectors determined directly from the $C^0$ interface representation, (10), pointwise expressions for the jump conditions for pressure (3) and and shear stress (4) are discontinuous between elements. To obtain a continuous representation of these jump conditions along the discretized boundary geometry, previous work by Kolahdouz et al.\cite{kolahdouz_immersed_2020} projected both pressure and shear stress jump conditions onto the same finite element basis used to represent the geometry. For example, the pressure jump condition in the continuous Lagrangian basis is
\begin{equation}
\llbracket  p\rrbracket_h^\text{C}(\mathbf{X}, t)  = \sum_{j = 1}^M \llbracket  p\rrbracket_j^\text{C}(t) \,
  \psi_j (\mathbf{X})
\end{equation}
in which $\llbracket  p\rrbracket_j^\text{C}(t)$ is the jump condition at the interface nodes at time $t$.
Projecting the pressure jump condition requires $\llbracket  p\rrbracket_h^\text{C}(\mathbf{X}, t)$ to satisfy
\begin{equation}
\int_{\Gamma_0} \llbracket  p\rrbracket_h^\text{C}(\mathbf{X}, t) \, \psi_j(\mathbf{X}) \, {\mathrm d} A = \int_{\Gamma_0} \jmath^{- 1} (\mathbf{X}, t)\, \mathbf{F}_h (\mathbf{X},
  t) \cdot \mathbf{n} (\mathbf{X}, t) \, \psi_j(\mathbf{X}) {\mathrm d} A , \forall j = 1, \ldots, M.
\end{equation}
If the geometry of the immersed interface is non-smooth, then $\mathbf{F}
(\mathbf{X}, t) \cdot \mathbf{n} (\mathbf{X}, t)$
is not well approximated by continuous piecewise linear basis functions. For such geometries, this projection produces errors that persist under grid refinement at sharp features in the true interface geometry, because at such features, the continuous jump condition $\mathbf{F}
(\mathbf{X}, t) \cdot \mathbf{n} (\mathbf{X}, t)$ is itself discontinuous.

If it is possible to align the mesh elements with the sharp features of the interface geometry, then the discontinuities in the jump condition will only occur at element boundaries. In this case, we can use a discontinuous basis for computing jump conditions. 
Discontinuous basis functions provide a better approximation for discontinuous fields and avoid $O (1)$ errors at the points of discontinuity. In this case, we can instead use discontinuous Lagrange polynomials
$\{ \zeta_j (\mathbf{X}) \}_{j = 1}^K$. $\zeta_j$ is a piecewise linear
Lagrange polynomial that is discontinuous across elements at shared nodes,
allowing nodal values to be different for neighboring elements. In contrast to Equation (11), we compute the projection into the discontinuous basis by requiring $\llbracket  p\rrbracket_h^\text{D}(\mathbf{X}, t)$ to satisfy

\begin{equation}
\int_{\Gamma_0} \llbracket  p\rrbracket_h^\text{D}(\mathbf{X}, t)\, \zeta_j(\mathbf{X}) \, {\mathrm d} A = \int_{\Gamma_0} \jmath^{- 1} (\mathbf{X}, t) \mathbf{F}_h (\mathbf{X},
  t) \cdot \mathbf{n} (\mathbf{X}, t) \, \zeta_j(\mathbf{X}){\mathrm d} A , \forall j = 1, \ldots, K.
\end{equation}
The projection of $\llbracket \frac{\partial \mathbf{u} (\mathbf{\boldsymbol{\chi}} (\mathbf{X}, t), t)}{\partial \mathbf{x}_i} \rrbracket$ is handled similarly. In this work, we evaluate these integrals using a seventh order Gaussian quadrature scheme.

\subsection{Finite Difference Approximation}

We use a staggered-grid discretization of the incompressible Navier-Stokes equations. The discretization approximates 
pressure at cell centers, and velocity and forcing terms at the center of cell edges (in two spatial dimensions) or cell faces (in three spatial dimensions).{\cite{griffith_volume_2012,griffith_accurate_2009}} Our computations use an isotropic grid, such that $h=\Delta x = \Delta y = \Delta z$. We use second-order accurate finite
difference stencils, and define the discrete divergence of the velocity
$\mathbf{D} \cdot \mathbf{u}$ on cell centers, and the discrete
pressure gradient $\mathbf{G} p$ and discrete Laplacian $L
\mathbf{u} $ on cell edges (or, in three spatial dimensions, faces). All numerical experiments use adaptive mesh refinement{\cite{griffith_immersed_2012}} to
increase computational efficiency. Finite difference stencils that cut through the interface require additional forcing terms, which are commonly referred to as correction terms, to account for discontinuities in the pressure and velocity gradients, as in previous work\cite{kolahdouz_immersed_2020}.

\subsection{Force Spreading and Velocity Interpolation}
The spreading operator, $\mathcal{S}[\boldsymbol{\chi} (\mathbf{X}, t)]$, relates the Lagrangian force
$\mathbf{F}_h (\mathbf{X}, t)  = \sum_{j = 1}^M  \mathbf{F}_j(t) \,\psi_j (\mathbf{X})$
to the Eulerian forces $\mathbf{f}$ on the Cartesian grid by summing all correction terms in the modified finite difference stencils involving the Lagrangian force. Correction terms, which are detailed in previous work \cite{kolahdouz_immersed_2020} depend on the projected jump conditions, as computed in Section 3.2.

The configuation of the Lagrangian mesh is updated by applying the no slip condition $\mathbf{U}_h = \mathcal{I}[{\boldsymbol{\chi},\mathbf{F}} ]\, \mathbf{u}$, in which $\mathcal{I}$ is a modified bilinear  (or trilinear in three spatial dimensions) interpolation in two spatial dimesions and subsequent
$L^2 $ projection of $\mathbf{u}$ about the mesh nodes. First, Eulerian velocities are interpolated to quadrature points along the interface, as described by Kolahdouz et al.\cite{kolahdouz_immersed_2020} These interpolated velocities are then used to compute the $L^2$ projection of the velocity at mesh nodes. These nodal velocities are represented in the finite element basis as $\mathbf{U}_h (\mathbf{X}, t)  = \sum_{j = 1}^M  \mathbf{U}_j(t) \,\psi_j (\mathbf{X})$. We use a CG projection for velocity interpolation because the velocity field is always continuous along the interface, even if the geometry includes sharp features.

\subsection{Time Integration}

Each step begins with known values of $\mathbf{\boldsymbol{\chi}}^n$ and $\mathbf{u}^n$ at time
$t^n$, and $p^{n - \frac{1}{2}}$ at time $t^{n - \frac{1}{2}}$. The goal is to
compute $\mathbf{\boldsymbol{\chi}}^{n + 1}, \mathbf{u}^{n + 1},$ and $p^{n +
\frac{1}{2}}$. First, an initial prediction of the structure location at time $t^{n + 1}$ is determined by
\begin{equation}
  \frac{\widehat{\mathbf{\boldsymbol{\chi}}}^{n + 1} - \mathbf{\boldsymbol{\chi}}^n}{\Delta t} =
  \mathbf{U}^n (\mathbf{\boldsymbol{\chi}}^n) = \mathcal{I}[\boldsymbol{\chi}^n,\mathbf{F}^n] \,\mathbf{u}^n  .
\end{equation}
We can then approximate the structure location at time $t^{n + \frac{1}{2}}$ by
\begin{equation}
  \mathbf{\boldsymbol{\chi}}^{n + \frac{1}{2}} = \frac{\widehat{\mathbf{\boldsymbol{\chi}}}^{n + 1}
  + \mathbf{\boldsymbol{\chi}}^n}{2}.
\end{equation}
Next, we solve for $\mathbf{u}^{n + 1}$ and $p^{n +
\frac{1}{2}}$ in which $\mathbf{f}^{n + \frac{1}{2}} = \mathcal{S}\,[\boldsymbol{\chi}^{n + \frac{1}{2}} ] \mathbf{F}(\boldsymbol{\chi}^{n + \frac{1}{2}},\mathbf{U}^n,t^{n + \frac{1}{2}})$:
\begin{eqnarray}
  \rho \left( \frac{\mathbf{u}^{n + 1} - \mathbf{u}^n}{\Delta t} +
  \mathbf{A}^{n+\frac12} \right) & = & - \mathbf{G} p^{n +
  \frac{1}{2}} + \mu L \left( \frac{\mathbf{u}^{n + 1} +
  \mathbf{u}^n}{2} \right) + \mathbf{f}^{n + \frac{1}{2}}, \\
  \mathbf{D} \cdot \mathbf{u}^{n + 1} & = & 0,
\end{eqnarray}
in which the non-linear advection term, $\mathbf{A}^{n + \frac{1}{2}} = \frac{3}{2}\mathbf{A}^{n} - \frac{1}{2}\mathbf{A}^{n-1}$, is handled with the xsPPM7 variant\cite{xsPPM7} of the piecewise
parabolic method.{\cite{colella_piecewise-parabolic_1982}} 
$\mathbf{G}, \mathbf{D} \cdot\mbox{}$, and $L$ are the discrete
gradient, divergence, and Laplacian operators. The system of equations is
iteratively solved via the FGMRES
algorithm with the projection method preconditioner\cite{griffith_accurate_2009}. Last, we update the structure's location, $\mathbf{\boldsymbol{\chi}}^{n + 1}$, via
\begin{equation}
  \frac{\mathbf{\boldsymbol{\chi}}^{n + 1} - \mathbf{\boldsymbol{\chi}}^n}{\Delta t}  = 
  \mathbf{U}^{n + \frac{1}{2}} =  \mathcal{I}[\boldsymbol{\chi}^{n+\frac{1}{2}},\mathbf{F}^{n+\frac{1}{2}}] \left(\frac{\mathbf{u}^n+ \mathbf{u}^{n+1}}{2}\right)  .
\end{equation}

\subsection{Software Implementation}

All computations for were completed through IBAMR \cite{ibamr}, which
utilizes parallel computing libraries and adaptive mesh refinement (AMR).
IBAMR uses other libraries for setting up meshes, fast linear algebra solvers,
and postprocessing, including SAMRAI \cite{samrai}, PETSc \cite{petsc,petsc2,petsc3}, \textit{hypre} \cite{hypre,hypre2}, and libMesh \cite{libmesh,libmesh2}.

\section{Numerical Experiments}

We first use flow past a stationary, rigid circular
cylinder as a benchmark case to examine the accuracy for smooth geometries, for which the original methodology produces accurate results. We then extend our geometries to a square
cylinder, and a wedge formed as a union of a square and rear-facing triangle. We also perform studies in three spatial dimensions, using a sphere and then a cube. All numerical examples are
non-dimensionalized.

We impose rigid body motions using a penalty method described in Section 2. Because the penalty parameter $\kappa$ is finite, generally there are discrepancies
between the prescribed and actual positions of the interface. In numerical tests, we choose $\kappa$ to ensure that the discrepancy in interface configurations $\varepsilon_{\mathbf{X}} = \| \mathbf{\xi} (\mathbf{X}, t) -
\mathbf{\boldsymbol{\chi}} (\mathbf{X}, t) \|$ satisfies $\varepsilon_{\mathbf{X}} \leq \frac{h}{4}$. For most experiments, we use block-structured
adaptive mesh refinement \cite{griffith_immersed_2012} with composite grids that include a total of 6 to 8 grid levels, with a refinement ratio of 2 between levels.
For all experiments, we set the Lagrangian mesh width to be twice as coarse as the background Cartesian grid spacing. 

\

\subsection{Flow Past a Circular Cylinder}
This section considers flow past a stationary cylinder of diameter $D=1$ centered at the origin of the computational domain $\Omega = [- 15, 45] \times [-30, 30]$ with $L=60$. 
For $\mathbf{u} = (u, v)$, we use incoming flow velocity $\mathbf{} u = 1$, and $v = \cos \left( \pi \frac{y}{L} \right) e^{- 2 t}$
to ensure that vortex shedding occurs at a consistent time. We set $\rho = 1$ and $\mu = \frac{1}{{Re}}$, in which $\text{Re}$ is the Reynolds number. We use
$\text{Re} = \frac{\rho U D}{\mu}=200$ to ensure that we observe vortex
shedding.  On the coarsest level
of the hierarchical grid, $n_\cc = 32$ grid cells are used in both directions.
The effective number of grid cells on the $\ell^{\text{th}}$ level is
${n_{\ell}}  = 2^{\ell - 1} n_\cc$. For convergence studies, we fix $n_\cc$ and vary $\ell_{\text{max}}$, the maximum number of grid levels, from six to eight. For comparisons to prior works, we use eight levels. The time step size $\Delta t$ scales with $h_{\text{finest}}$ such that $\Delta t = \frac{
h_{\text{finest}}}{20}$, in which $h_{\text{finest}} = \frac{L}{{n_{\ell {\text{max}}}}}$. For this time step size and grid spacing, the Courant-Friedrichs-Lewy (CFL) number is approximately 0.04 once the model reaches periodic steady state. $\kappa$ and $\eta$ are scaled with the grid resolution on the finest level, such that $\kappa = \frac{C_\kappa}{h_{\text{finest}}} $ and $\eta = \frac{C_\eta}{h_{\text{finest}}}$ with $C_\kappa = 3.413$ and $ C_\eta = 0.025 $ using a bisection method. The outflow
boundary uses zero normal and tangential traction, and the top and bottom boundaries use zero tangential traction and $v = 0$.
Figure 2 shows a snapshot of the vorticity field, $\omega$, at $t = 100$ for both CG
and DG jump condition handling. Both flow fields show the expected vortex shedding and a max$(\varepsilon_{\mathbf{X}})$ on the same order of magnitude.

\begin{figure}[h]
\center{}
\resizebox{300pt}{!}{\includegraphics[]{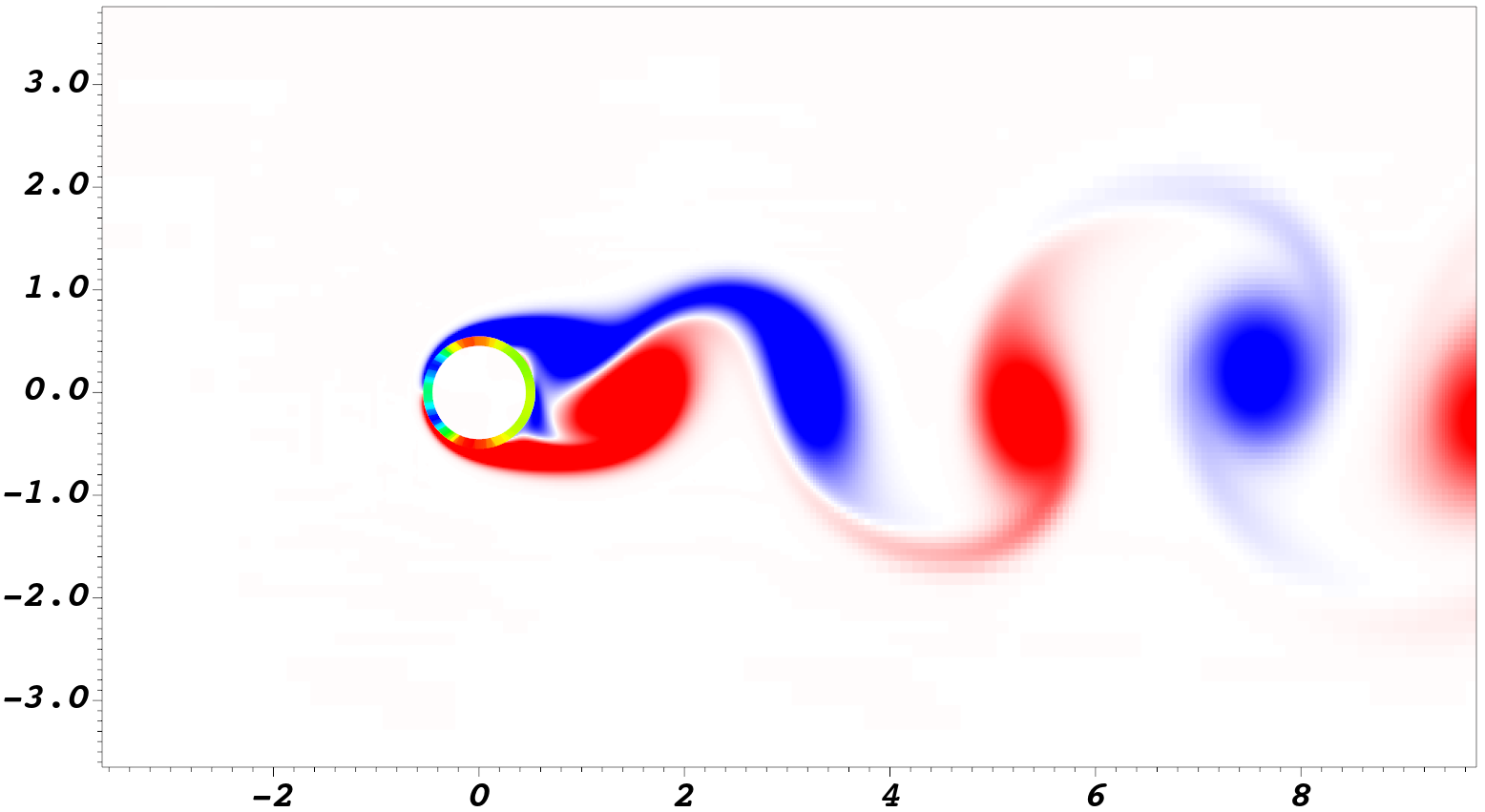}}
\resizebox{70pt}{60pt}{
\tikzset{every picture/.style={line width=0.75pt}} 

\begin{tikzpicture}[x=0.75pt,y=0.75pt,yscale=-1,xscale=1]

\draw (78.97,110.76) node  {\includegraphics[width=50.96pt,height=79.62pt]{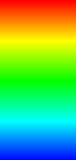}};
\draw [color={rgb, 255:red, 0; green, 0; blue, 0 }  ,draw opacity=1 ]   (112.94,59.2) -- (136.64,59.2) ;
\draw    (112.94,163.84) -- (136.64,163.84) ;
\draw    (112.94,110.76) -- (136.64,110.76) ;
\draw [color={rgb, 255:red, 0; green, 0; blue, 0 }  ,draw opacity=1 ]   (112.94,84.98) -- (136.64,84.98) ;
\draw    (112.94,136.54) -- (136.64,136.54) ;

\draw (70,40) node [anchor=north west][inner sep=0.75pt]  [font=\Large]  {$\varepsilon $};
\draw (147.82,45.42) node [anchor=north west][inner sep=0.75pt]   [align=left] {3.6e-3};
\draw (147.82,72.72) node [anchor=north west][inner sep=0.75pt]   [align=left] {2.8e-3};
\draw (147.82,100.02) node [anchor=north west][inner sep=0.75pt]   [align=left] {2.1e-3};
\draw (147.82,125.8) node [anchor=north west][inner sep=0.75pt]   [align=left] {1.3e-3};
\draw (147.82,151.58) node [anchor=north west][inner sep=0.75pt]   [align=left] {0.000};

\end{tikzpicture}}

\center{}
\resizebox{300pt}{!}{\includegraphics[]{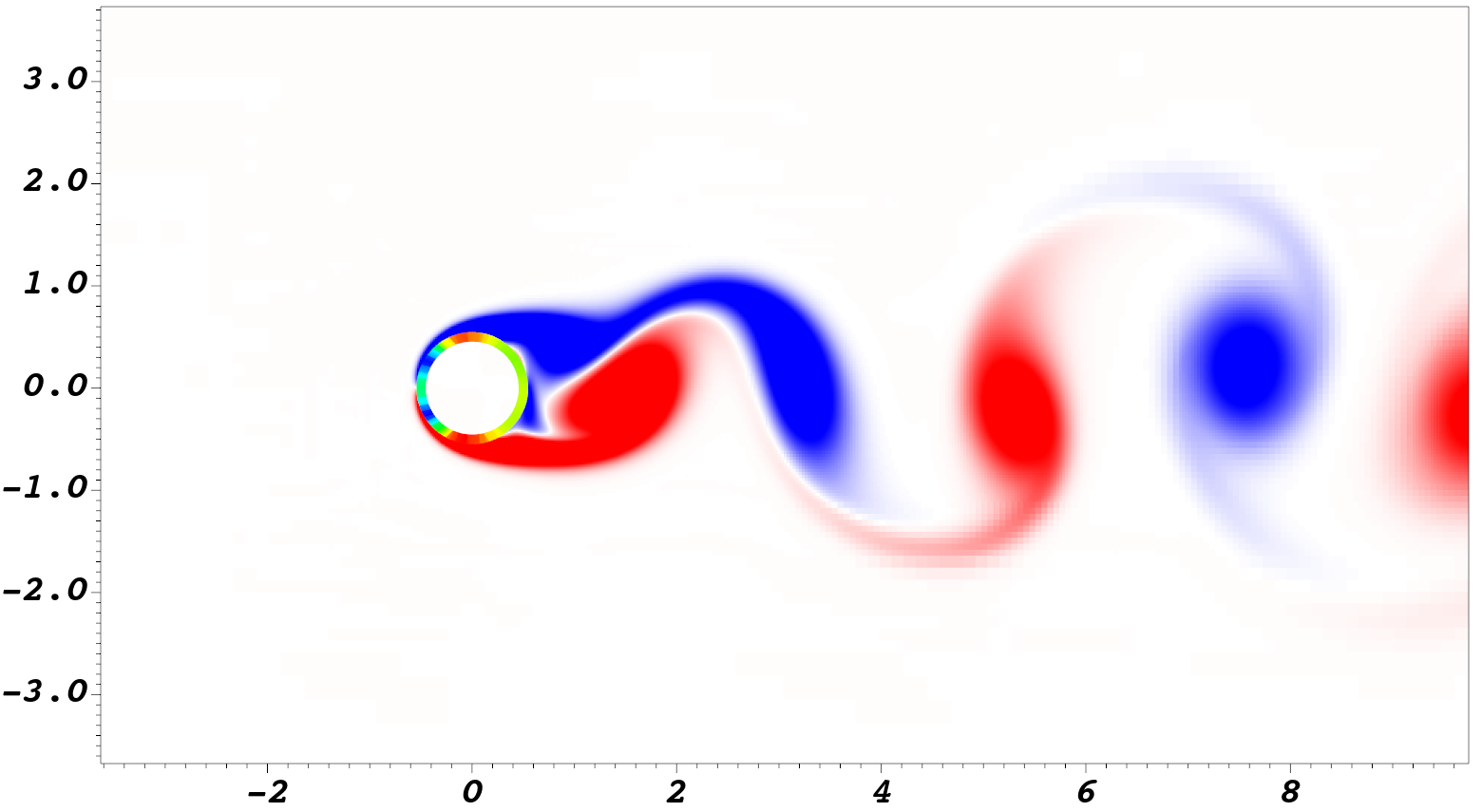}}
\resizebox{70pt}{120pt}{
\tikzset{every picture/.style={line width=0.75pt}} 

\begin{tikzpicture}[x=0.75pt,y=0.75pt,yscale=-1,xscale=1]

\draw (83.97,256.76) node  {\includegraphics[width=50.96pt,height=79.62pt]{Figures/colorbar_displacement.png}};
\draw [color={rgb, 255:red, 0; green, 0; blue, 0 }  ,draw opacity=1 ]   (117.94,205.2) -- (141.64,205.2) ;
\draw    (117.94,309.84) -- (141.64,309.84) ;
\draw    (117.94,256.76) -- (141.64,256.76) ;
\draw [color={rgb, 255:red, 0; green, 0; blue, 0 }  ,draw opacity=1 ]   (117.94,230.98) -- (141.64,230.98) ;
\draw    (117.94,282.54) -- (141.64,282.54) ;

\draw [color={rgb, 255:red, 0; green, 0; blue, 0 }  ,draw opacity=1 ]   (116.94,58.62) -- (140.64,58.62) ;
\draw    (116.94,163.26) -- (140.64,163.26) ;
\draw    (116.94,110.18) -- (140.64,110.18) ;
\draw [color={rgb, 255:red, 0; green, 0; blue, 0 }  ,draw opacity=1 ]   (116.94,84.4) -- (140.64,84.4) ;
\draw    (116.94,135.96) -- (140.64,135.96) ;
\draw (82,111.05) node  {\includegraphics[width=52.5pt,height=80.92pt]{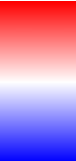}};

\draw (151.82,151) node [anchor=north west][inner sep=0.75pt]   [align=left] {\mbox{-}2.000};
\draw (151.82,125.22) node [anchor=north west][inner sep=0.75pt]   [align=left] {\mbox{-}1.000};
\draw (151.82,99.44) node [anchor=north west][inner sep=0.75pt]   [align=left] {0.000};
\draw (151.82,72.14) node [anchor=north west][inner sep=0.75pt]   [align=left] {1.000};
\draw (151.82,44.85) node [anchor=north west][inner sep=0.75pt]   [align=left] {2.000};
\draw (71,40) node [anchor=north west][inner sep=0.75pt]  [font=\large]  {$\omega $};
\draw (152.82,297.58) node [anchor=north west][inner sep=0.75pt]   [align=left] {0.000};
\draw (152.82,271.8) node [anchor=north west][inner sep=0.75pt]   [align=left] {1.3e-3};
\draw (152.82,246.02) node [anchor=north west][inner sep=0.75pt]   [align=left] {2.1e-3};
\draw (152.82,218.72) node [anchor=north west][inner sep=0.75pt]   [align=left] {2.8e-3};
\draw (152.82,191.42) node [anchor=north west][inner sep=0.75pt]   [align=left] {3.6e-3};
\draw (75,184) node [anchor=north west][inner sep=0.75pt]  [font=\Large]  {$\varepsilon $};

\end{tikzpicture}
  }
  \
  \caption{The vorticity field and the variable \(\varepsilon_{\mathbf{X}}\) are visualized over a subset of the fluid domain at \(t = 100\). The top panel uses DG-IIM, while the bottom panel uses CG-IIM jump conditions.}
\end{figure}

\begin{figure}[h]

  \resizebox{220pt}{!}{\includegraphics{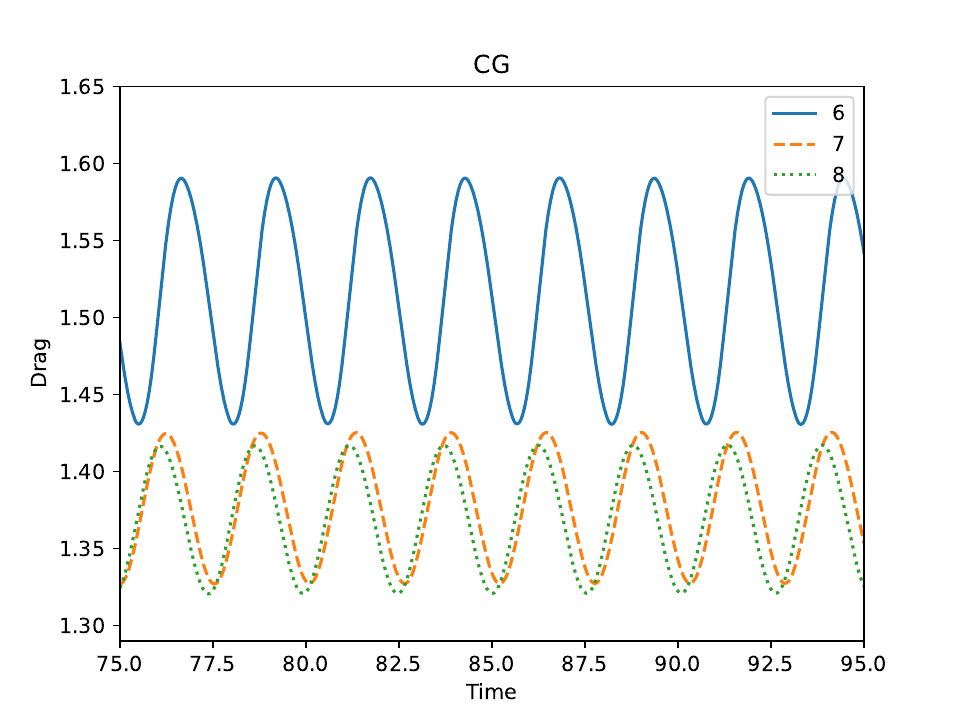}}\resizebox{220pt}{!}{\includegraphics{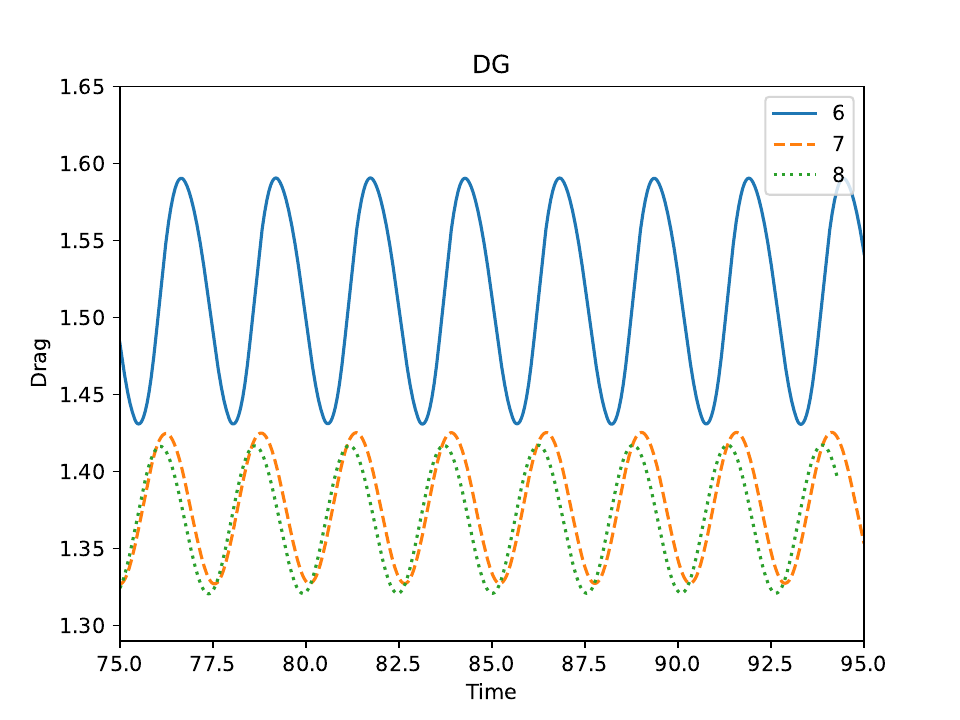}}

    \resizebox{220pt}{!}{\includegraphics{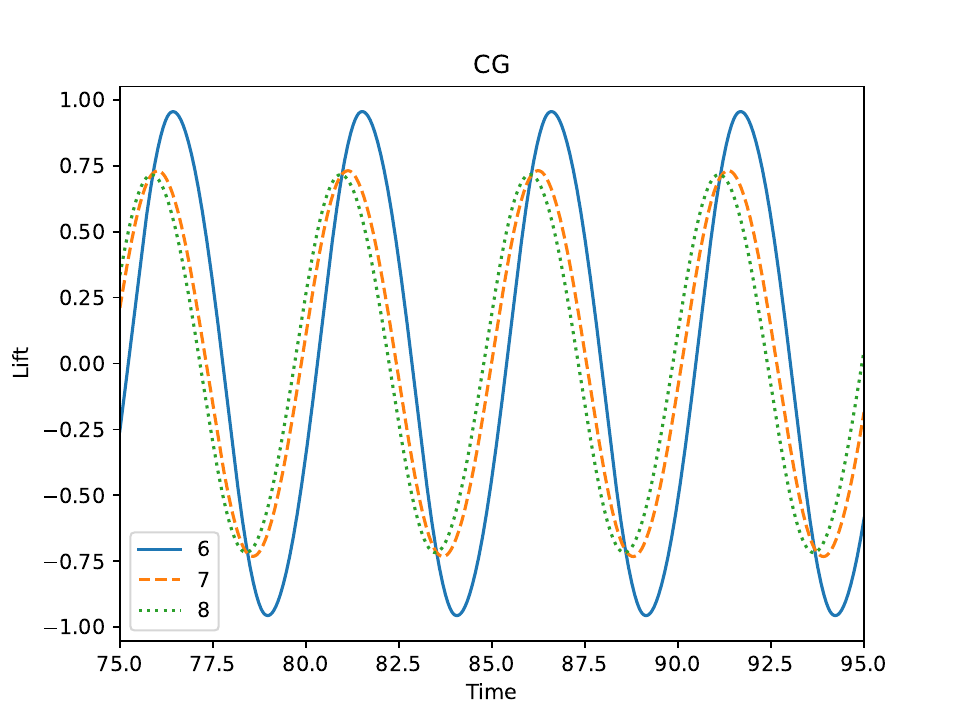}}\resizebox{220pt}{!}{\includegraphics{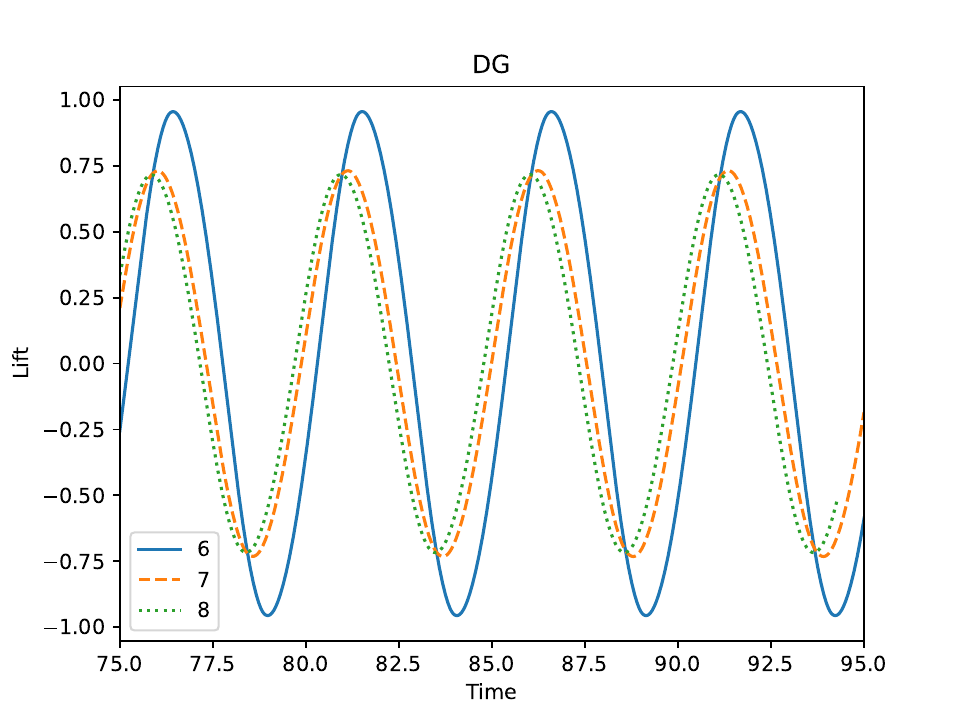}}
  \caption{Time histories of lift and drag for flow past a circular cylinder for $\ell_{\text{max}} = \{6,7,8\}$. Both CG-IIM and DG-IIM exhibit convergence under grid refinement and exhibit similar accuracy at comparable grid resolutions.}
\end{figure}

We quantify the simulation outputs using the nondimensionalized drag and lift coefficients, which are computed via
\begin{equation}
  C_{\text{D}} = \frac{- \int_{\Gamma_0} F^x (\mathbf{X}, t)
  \, {\mathrm d}s}{\frac{1}{2} \rho U^2 D}
\end{equation}
\begin{equation}
  {C_{\text{L}}}_{\textrm{}} = \frac{- \int_{\Gamma_0} F^y (\mathbf{X}, t)
  \, {\mathrm d}s}{\frac{1}{2} \rho U^2 D}
\end{equation}
As a verification of the method, we conduct a grid convergence test (Figure 3). For both lift and drag coefficients, DG-IIM and CG-IIM both clearly converge to quantities within the range of previous literature values. 
For a discrete
collection of times $\{ t_i \}_{i = 1}^M$ and the corresponding drag or lift
coefficients and $C (t_i)$, we compute the average coefficient as $C_{\text{avg}} = \frac{\max(C)+\min(C)}{2}$. Root-mean-squared drag and lift are computed as

\begin{equation}
C_{\text{rms}} = \sqrt{\frac{1}{M} \sum_{i = 1}^M \left(
C (t_i) - C_{\text{avg}} \right)^2}.
\end{equation}
The Strouhal number (St) is computed as $\text{St} = \frac{f D }{U}$, in which $f$ is the frequency of vortex shedding, $D$ is the diameter of the cylinder, and $U$ is the horizontal inflow velocity.

\begin{figure}[H]
  \resizebox{230pt}{!}{\includegraphics{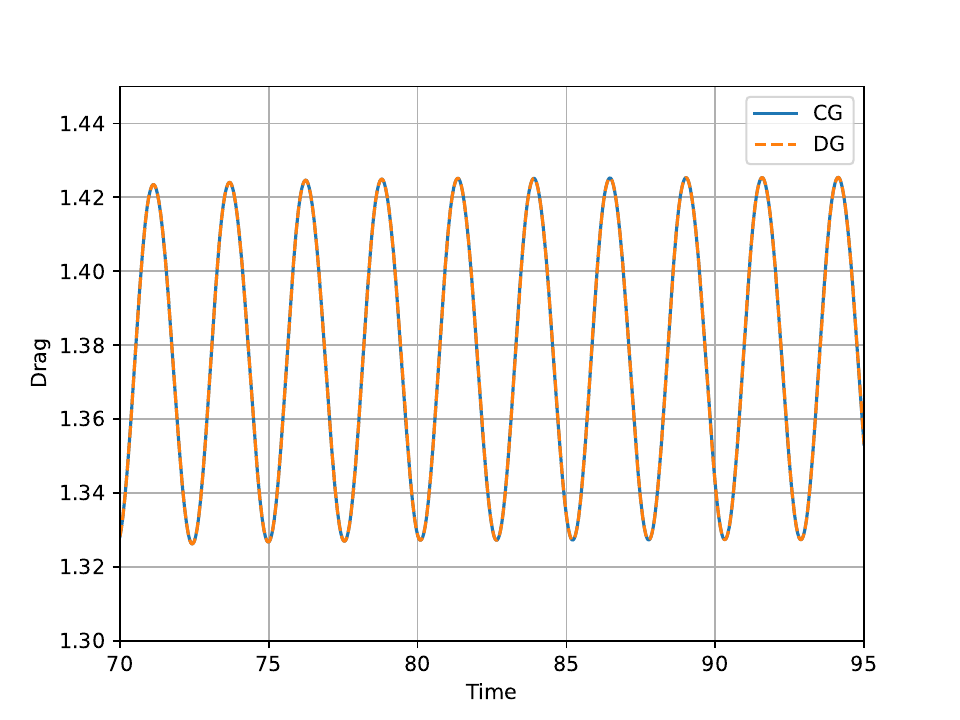}}\resizebox{230pt}{!}{\includegraphics{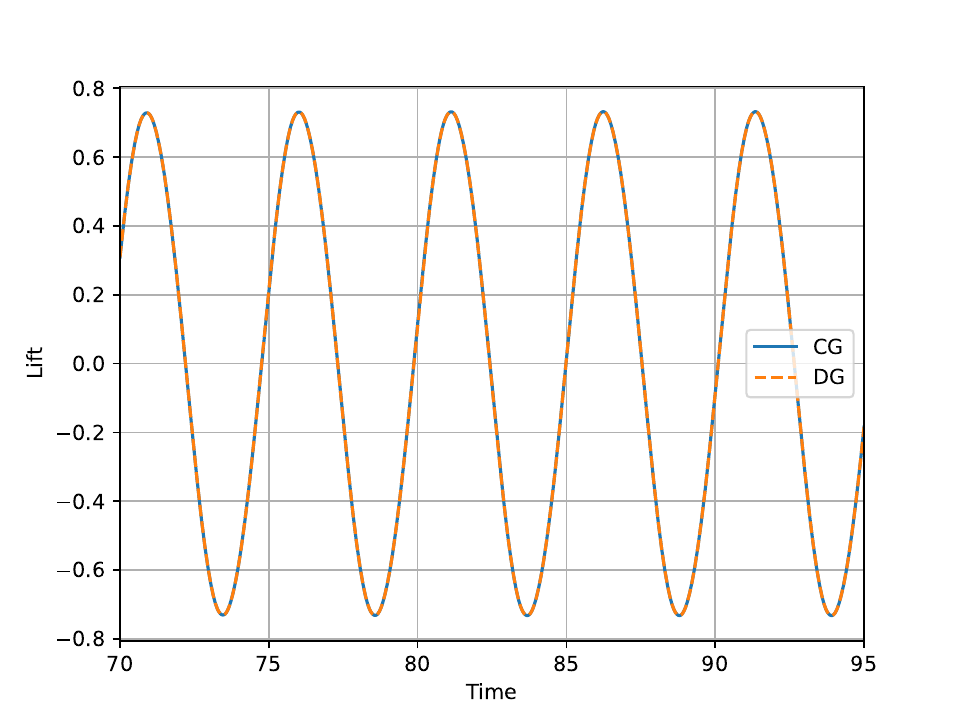}}
  \caption{Comparison of drag (left) and lift (right) coefficients over time
  for both CG and DG formulations for circular cylinder flow with 7 levels of
  refinement. The absolute pointwise difference in  $C_{\Davg}$, $C_{\Drms}$, $C_{\Lmax}$, and $C_{\Lrms}$
  between the two formulations are all less than 3e-5. }
\end{figure}
Figure 4 shows an agreement of these integrated quantities between the two methods. 
\begin{table}[t]
\center{}
  \begin{tabular}{l}
    \begin{tabular}{|c|c|c|c|}
      \hline
      \textbf{Reference} & $C_\text{D}$ & $C_\text{L}$ & St\\
      \hline
      Griffith and Luo{\cite{griffith_hybrid_2017}} & $1.360 \pm 0.046$ & $\pm
      0.70$ & 0.195\\
      \hline
      Liu et al. {\cite{liu_preconditioned_1998}} & $1.310 \pm 0.049$ & $\pm
      0.69$ & 0.192\\
      \hline
      Braza et al. {\cite{braza_numerical_1986}} & $1.400 \pm 0.050$ & $\pm
      0.75$ & 0.200\\
      \hline
      Calhoun{\cite{calhoun_cartesian_2002}} & $1.172 \pm 0.058$ & $\pm
      0.67$ & 0.202\\
      \hline
      Present CG & $1.370 \pm 0.049$ &
      $\pm 0.72$ & 0.197\\
      \hline
      Present DG & $1.370 \pm 0.049$ & $\pm 0.72$& 0.197\\
      \hline
    \end{tabular} 
  \end{tabular}
  \caption{Comparisons of lift and drag coefficients for $\text{Re} = 200$ flow
  past a circular cylinder to results from prior studies.  }
\end{table}
The present lift and drag values fall within the
range of literature values (Table 1). The study by Griffith and Luo used an immersed
boundary method with finite elements for structure and finite difference
approximations for Eulerian variables. The study by Liu et al. used a
finite differences with a turbulence model. Braza et al. used finite
differences on the Navier Stokes equation specifically for cylinder flow.
Lastly, the study by Calhoun used an IIM with finite differences on a Cartesian grid{\cite{calhoun_cartesian_2002}}.

\subsection{Flow Past a Square Cylinder}
This section considers flow past a stationary square cylinder of diameter $D=1$ centered at the origin of the computational domain $\Omega = [- 15, 45] \times [-30, 30]$ with $L=60$. 
For $\mathbf{u} = (u, v)$, we use incoming flow velocity $\mathbf{} u = 1$, and $v = \cos \left( \pi \frac{y}{L} \right) e^{- 2 t}$
to ensure that vortex shedding occurs at a consistent time. We set $\rho = 1$ and $\mu = \frac{1}{\text{Re}}$.  We use
$\text{Re} = \frac{\rho U D}{\mu}=100$ to ensure that we observe vortex
shedding. 
The effective number of grid cells on the $\ell^{\text{th}}$ level is
${n_{\ell}}  = 2^{\ell - 1} n_\cc$,  with $\ell_{\text{max}}=6$. For convergence studies, we vary $n_\cc = 64$, $128$, $256$, and $512$. Comparisons to prior works, use the finest resolution considered, $n_\cc = 512$. The time step size $\Delta t$ scales with $h_{\text{finest}}$ such that $\Delta t = \frac{
h_{\text{finest}}}{20}$, in which $h_{\text{finest}} = \frac{L}{{n_{\ell {\text{max}}}}}$. For this time step size and grid spacing, the CFL number is approximately 0.05 once the model reaches periodic steady state. $\kappa$ and $\eta$ are scaled with the grid resolution on the finest level, such that $\kappa = \frac{C_\kappa}{h_{\text{finest}}} $ and  $\eta = \frac{C_\eta}{h_{\text{finest}}}$. We set $C_\kappa = 3.413$ and $ C_\eta = 0.0036$, which were empirically determined using bisection. The outflow
boundary uses zero normal and tangential traction, and the top and bottom boundaries use zero tangential traction and $v = 0$.
Figure 5 gives a snapshot of the vorticity field and interface at $t = 150$.

\begin{figure}[H]
\center{}
  \raisebox{0.0\height}{\includegraphics[width=10.4553325462416cm,height=5.80850714941624cm]{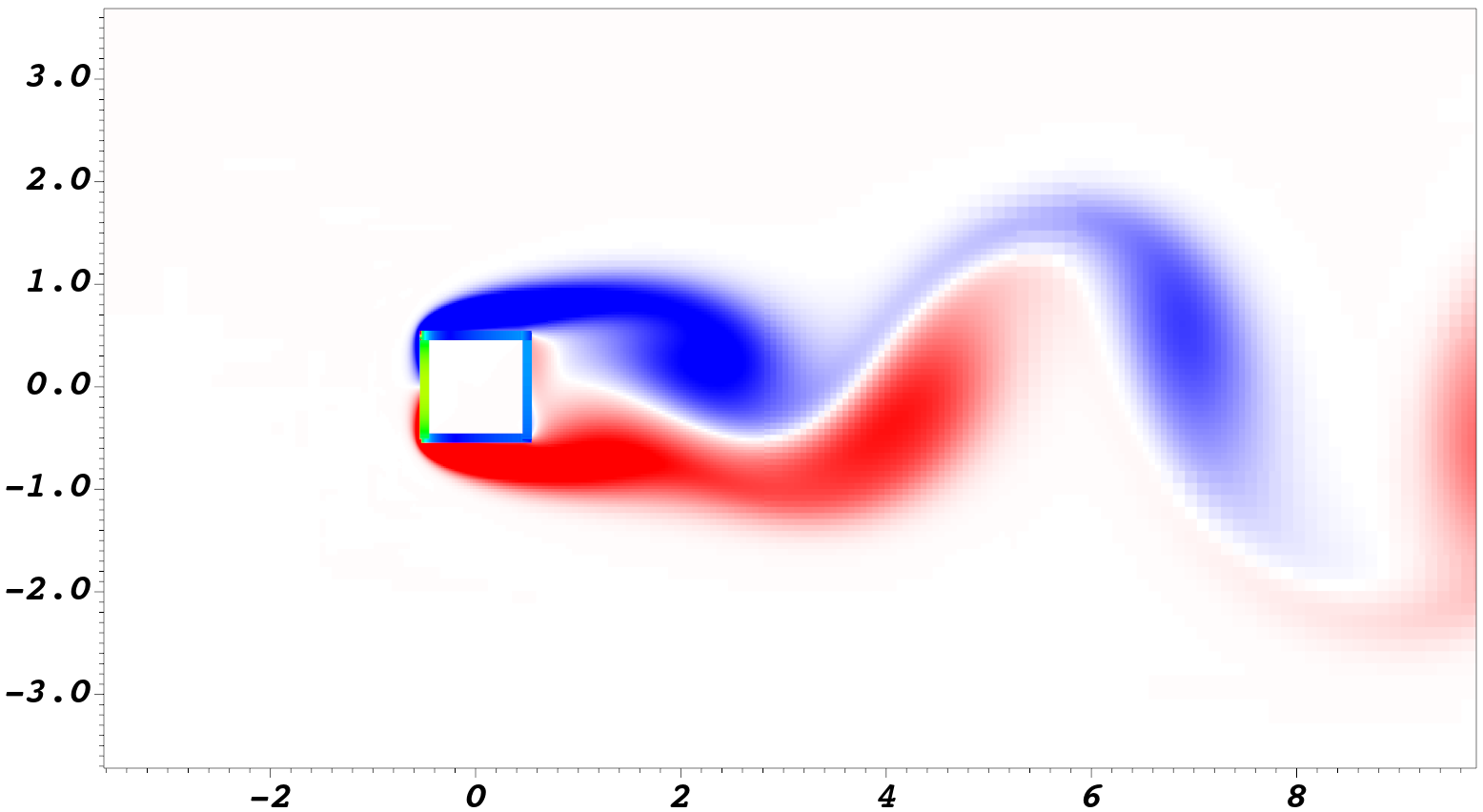}}\resizebox{80pt}{!}{

\tikzset{every picture/.style={line width=0.75pt}} 

\begin{tikzpicture}[x=0.75pt,y=0.75pt,yscale=-1,xscale=1]

\draw (78.97,110.76) node  {\includegraphics[width=50.96pt,height=79.62pt]{Figures/colorbar_displacement.png}};
\draw [color={rgb, 255:red, 0; green, 0; blue, 0 }  ,draw opacity=1 ]   (112.94,59.2) -- (136.64,59.2) ;
\draw    (112.94,163.84) -- (136.64,163.84) ;
\draw    (112.94,110.76) -- (136.64,110.76) ;
\draw [color={rgb, 255:red, 0; green, 0; blue, 0 }  ,draw opacity=1 ]   (112.94,84.98) -- (136.64,84.98) ;
\draw    (112.94,136.54) -- (136.64,136.54) ;

\draw (70,40) node [anchor=north west][inner sep=0.75pt]  [font=\Large]  {$\varepsilon $};
\draw (147.82,45.42) node [anchor=north west][inner sep=0.75pt]   [align=left] {3.6e-3};
\draw (147.82,72.72) node [anchor=north west][inner sep=0.75pt]   [align=left] {2.8e-3};
\draw (147.82,100.02) node [anchor=north west][inner sep=0.75pt]   [align=left] {2.1e-3};
\draw (147.82,125.8) node [anchor=north west][inner sep=0.75pt]   [align=left] {1.3e-3};
\draw (147.82,151.58) node [anchor=north west][inner sep=0.75pt]   [align=left] {0.000};

\end{tikzpicture}}

  \center{}
  \raisebox{0.0\height}{\includegraphics[width=10.4553325462416cm,height=5.84920634920635cm]{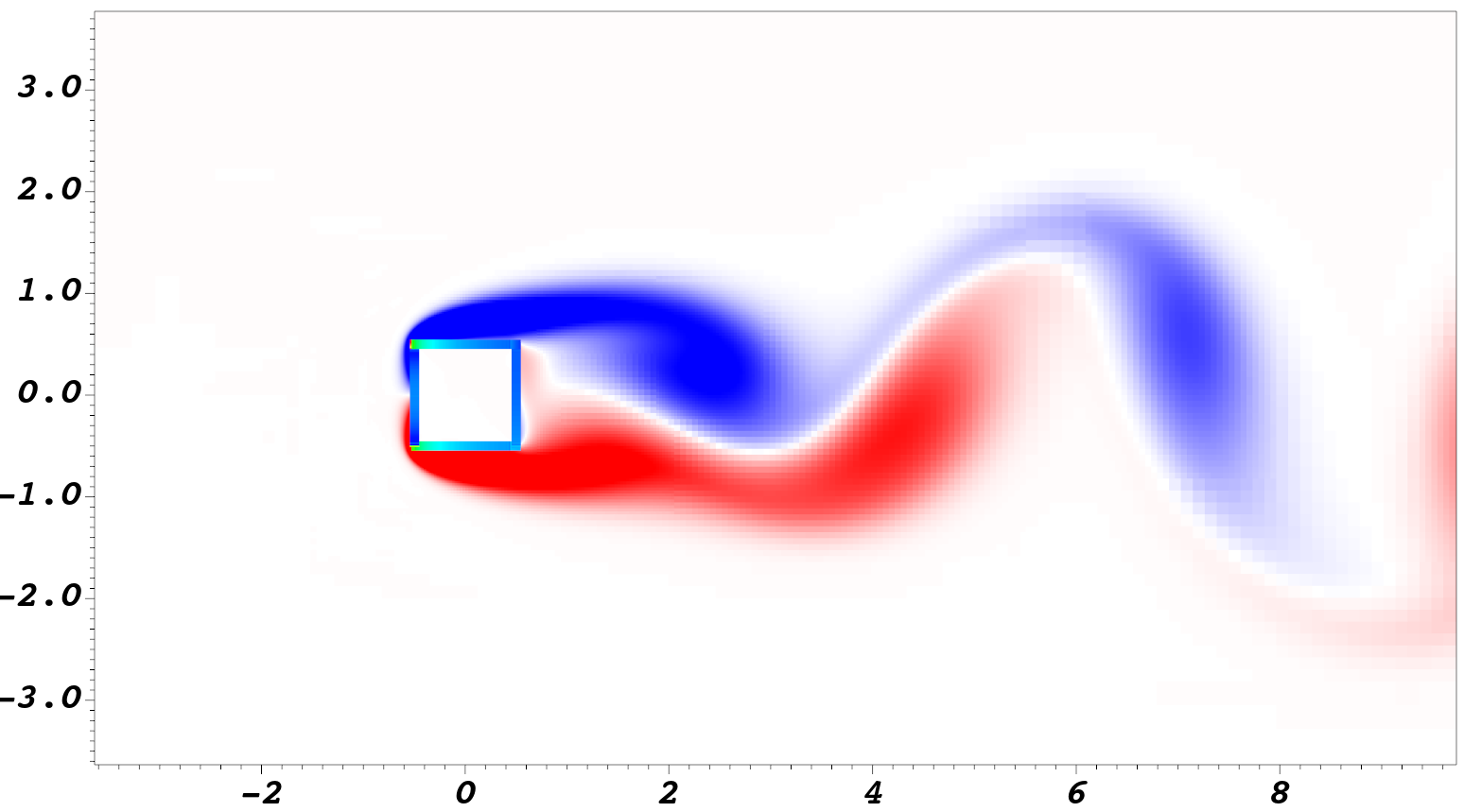}}\resizebox{80pt}{!}{

\tikzset{every picture/.style={line width=0.75pt}} 

\begin{tikzpicture}[x=0.75pt,y=0.75pt,yscale=-1,xscale=1]

\draw (83.97,256.76) node  {\includegraphics[width=50.96pt,height=79.62pt]{Figures/colorbar_displacement.png}};
\draw [color={rgb, 255:red, 0; green, 0; blue, 0 }  ,draw opacity=1 ]   (117.94,205.2) -- (141.64,205.2) ;
\draw    (117.94,309.84) -- (141.64,309.84) ;
\draw    (117.94,256.76) -- (141.64,256.76) ;
\draw [color={rgb, 255:red, 0; green, 0; blue, 0 }  ,draw opacity=1 ]   (117.94,230.98) -- (141.64,230.98) ;
\draw    (117.94,282.54) -- (141.64,282.54) ;

\draw [color={rgb, 255:red, 0; green, 0; blue, 0 }  ,draw opacity=1 ]   (116.94,58.62) -- (140.64,58.62) ;
\draw    (116.94,163.26) -- (140.64,163.26) ;
\draw    (116.94,110.18) -- (140.64,110.18) ;
\draw [color={rgb, 255:red, 0; green, 0; blue, 0 }  ,draw opacity=1 ]   (116.94,84.4) -- (140.64,84.4) ;
\draw    (116.94,135.96) -- (140.64,135.96) ;
\draw (82,111.05) node  {\includegraphics[width=52.5pt,height=80.92pt]{Figures/vorticity_colorbar.png}};

\draw (151.82,151) node [anchor=north west][inner sep=0.75pt]   [align=left] {\mbox{-}2.000};
\draw (151.82,125.22) node [anchor=north west][inner sep=0.75pt]   [align=left] {\mbox{-}1.000};
\draw (151.82,99.44) node [anchor=north west][inner sep=0.75pt]   [align=left] {0.000};
\draw (151.82,72.14) node [anchor=north west][inner sep=0.75pt]   [align=left] {1.000};
\draw (151.82,44.85) node [anchor=north west][inner sep=0.75pt]   [align=left] {2.000};
\draw (71,40) node [anchor=north west][inner sep=0.75pt]  [font=\large]  {$\omega $};
\draw (152.82,297.58) node [anchor=north west][inner sep=0.75pt]   [align=left] {0.000};
\draw (152.82,271.8) node [anchor=north west][inner sep=0.75pt]   [align=left] {1.6e-3};
\draw (152.82,246.02) node [anchor=north west][inner sep=0.75pt]   [align=left] {3.1e-3};
\draw (152.82,218.72) node [anchor=north west][inner sep=0.75pt]   [align=left] {4.7e-3};
\draw (152.82,191.42) node [anchor=north west][inner sep=0.75pt]   [align=left] {6.3e-3};
\draw (75,184) node [anchor=north west][inner sep=0.75pt]  [font=\Large]  {$\varepsilon $};

\end{tikzpicture}}

  \
  \caption{The vorticity field and the variable \(\varepsilon_{\mathbf{X}}\) are visualized over a subset of the fluid domain at \(t = 150\). The top panel uses DG-IIM, while the bottom panel uses CG-IIM jump conditions.}
\end{figure}

Figure 6 compare the discrepancy in interface configurations between CG-IIM and DG-IIM for various grid resolutions. Across all resolutions, the steady-state $\max (\varepsilon_{\mathbf{X}})$ for CG-IIM is twice as large as DG-IIM.
\begin{figure}[H]
    \center{}
  \resizebox{300pt}{!}{\includegraphics{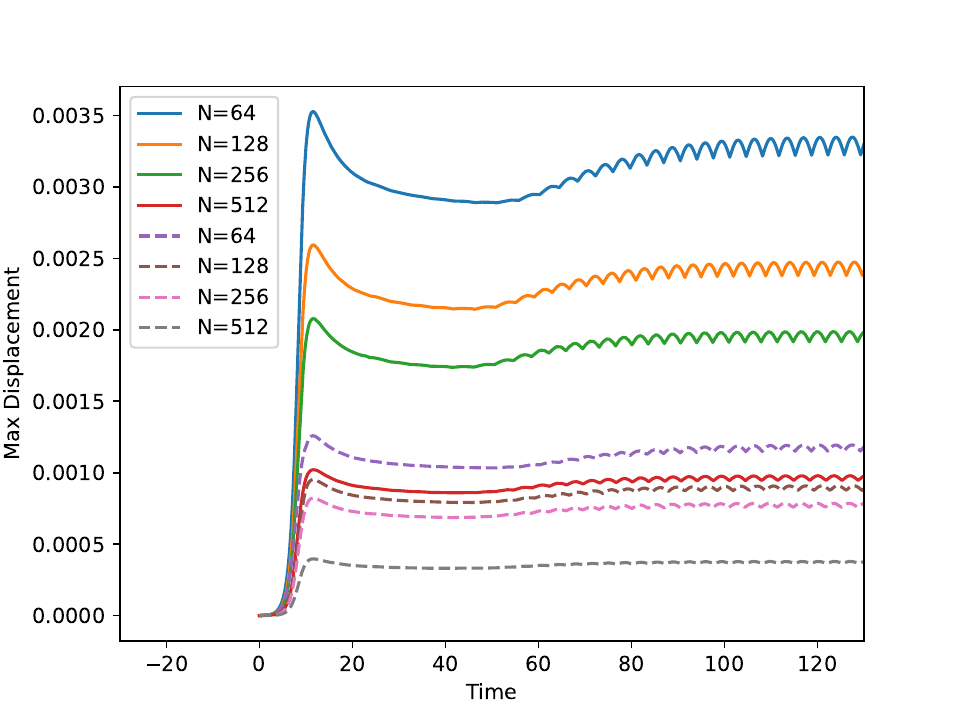}}
  
  \caption{Flow past a stationary square cylinder time series for $\max (\varepsilon_{\mathbf{X}})$ using DG-IIM and CG-IIM formulations with $\kappa = \frac{C_\kappa}{h_{\text{finest}}} $, $\eta = \frac{C_\eta}{h_{\text{finest}}}$, and $n_\cc = \{64,128,256,512\}$ grid cells on the coarsest level with six levels of refinement. CG-IIM results use solid lines, and the DG-IIM results use dashed lines. The CG-IIM formulation has larger displacement for all grid discretizations. }
\end{figure}
Next, as detailed in Figure 7, we study the convergence of the drag and lift coefficients under grid refinement. 

\begin{figure}[H]

  \resizebox{220pt}{!}{\includegraphics{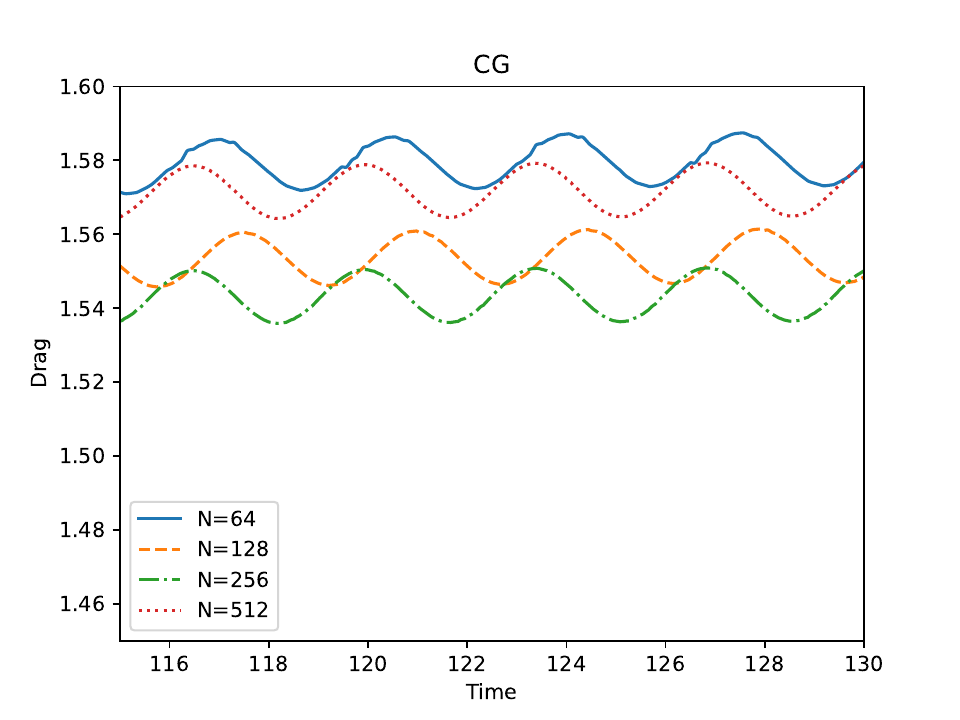}}\resizebox{220pt}{!}{\includegraphics{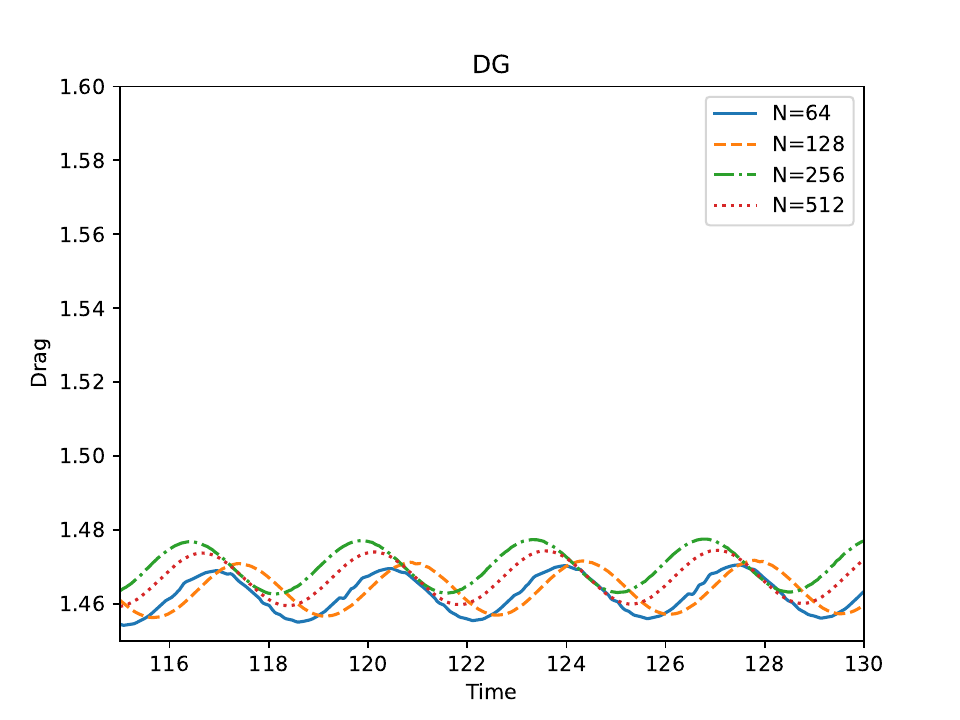}}
  \resizebox{220pt}{!}{\includegraphics{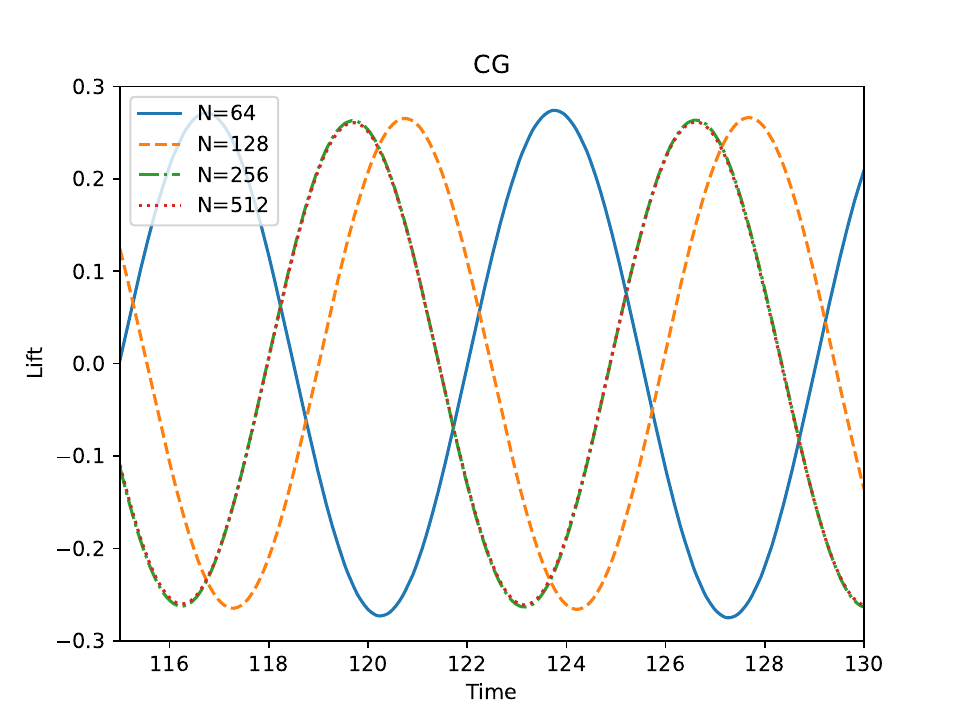}}\resizebox{220pt}{!}{\includegraphics{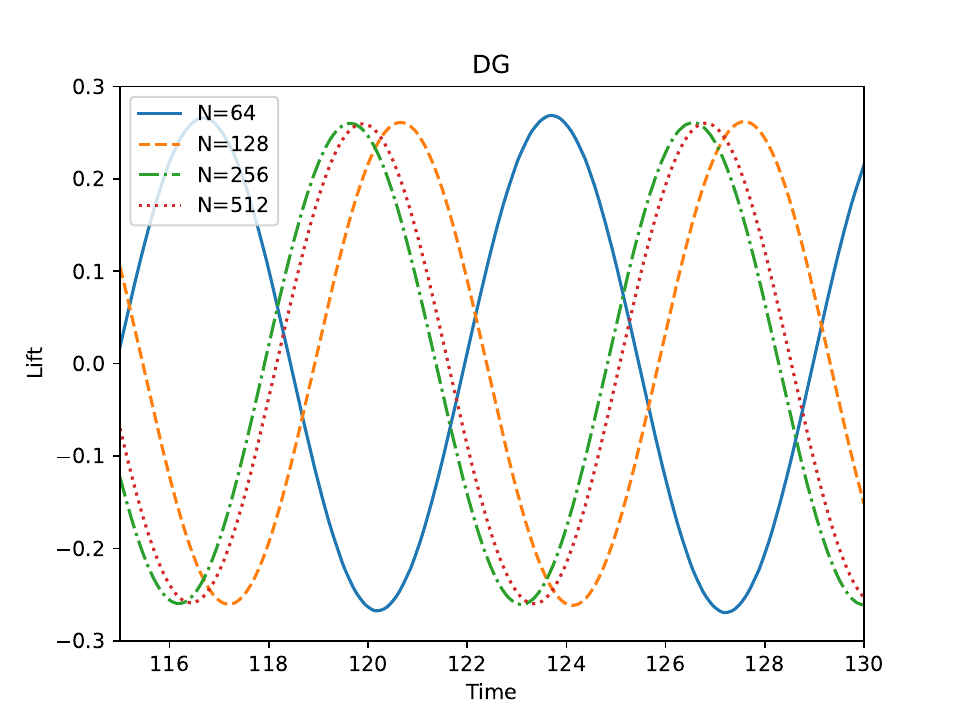}}

  \caption{Time histories of lift and drag coefficients for flow past a square cylinder with $n_\cc = \{64,128,256,512\}$. The lift coefficient exhibits convergence for both methods. In contrast, the drag coefficient converges to literature values only while using DG-IIM. }
\end{figure}

Both CG-IIM and DG-IIM agree for lift coefficients across all grid resolutions. The DG-IIM drag coefficient converges with grid refinement, whereas the CG-IIM formulation does not achieve convergence on practical grids. 
We compare our most refined periodic steady state
values to results from prior studies for $\text{Re}=100$ (Table 2). $C_{\text{L}_{\text{rms}}}$ agrees with prior studies for both formulations. $C_\Davg$ agrees with the range of literature values for DG-IIM, whereas CG-IIM yields a larger value. The Strouhal number for the DG-IIM formulation is closer to literature values than CG-IIM, that exhibits extremely slow convergence under grid refinement. 

\begin{table}[H]
\center{}
  \begin{tabular}{l}
    \begin{tabular}{|c|c|c|c|}
      \hline
      \textbf{Reference} & $C_\Davg$ & $C_{\text{L}_{\text{rms}}}$ & St\\
      \hline
      Sohankar et al.{\cite{sohankar_numerical_1997}} & 1.43 & $0.15$ &
      0.146\\
      \hline
      Cheng et al.{\cite{cheng_numerical_2007}} & $1.44$ & $0.15$ & 0.144\\
      \hline
      Lam et al.{\cite{lam_numerical_2012}} & 1.51 & 0.18 & 0.141\\
      \hline
      Berrone et al.{\cite{berrone_numerical_2011}} & 1.48 & 0.18 & 0.145\\
      \hline
      Sen et al.{\cite{sen_flow_2011}} & 1.52 & 0.19 & 0.145\\
      \hline
      Ryu \& Iaccarino{\cite{ryu_vortex-induced_2017}} & 1.56 & 0.19 &
      0.141\\
      \hline
      Present CG & 1.57 & 0.18 & 0.138\\
      \hline
      Present DG & 1.47 & 0.18 & 0.140\\
      \hline
    \end{tabular} 
  \end{tabular}
  \caption{Lift and drag coefficients, and Strouhal number comparison for
  $\text{Re} = 100$ flow past a square cylinder. }
\end{table}

Last, we compare our DG-IIM coefficients to values reported by Sen et al{\cite{sen_flow_2011}},  who performed numerical experiments for flow past a square cylinder at a variety of Reynolds numbers. Table 3 shows close
agreement in present calculations with those integrated quantities. Table 4 shows agreement in the
circulation lengths for steady flows with cylinder diameter $D = 1$.

\begin{table}[H]
\center{}
  \begin{tabular}{l}
    \begin{tabular}{|c|c|c|c|c|}
      \hline
      \text{Re} & \multicolumn{2}{c|}{$C_{\text{L}_{\max}}$}  & \multicolumn{2}{c|}{$C_{\text{L}_{\text{rms}}}$} \\
      \hline
      & Sen et al. & Present DG& Sen et al. & Present DG\\
      \hline
      60 & 0.13 & 0.124 & 0.09 & 0.0839\\
      \hline
      80 & 0.21 & 0.209 & 0.15 & 0.147\\
      \hline
      100 & 0.27 & 0.273 & 0.19 & 0.192\\
      \hline
      120 & 0.33 & 0.335 & 0.23 & 0.237\\
      \hline
      140 & 0.39 & 0.412 & 0.28 & 0.292\\
      \hline
    \end{tabular} \\
    \begin{tabular}{|c|c|c|c|c|}
      \hline
      \text{Re} & \multicolumn{2}{c|}{$C_{\text{D}_\text{avg}}$}  & \multicolumn{2}{c|}{$C_{\text{D}_{\text{rms}}}$} \\
      \hline
      & Sen et al. & Present DG & Sen et al. & Present DG\\
      \hline
      60 & 1.64 & 1.55 & 0.0009 & 0.0009\\
      \hline
      80 & 1.57 & 1.50 & 0.0028 & 0.0027\\
      \hline
      100 & 1.53 & 1.47 & 0.0055 & 0.0055\\
      \hline
      120 & 1.51 & 1.46 & 0.0096 & 0.0092\\
      \hline
      140 & 1.50 & 1.47 & 0.0143 & 0.0144\\
      \hline
    \end{tabular}
  \end{tabular}
  
  \
  \caption{Comparison of $C_\Lmax$, $C_\Lrms$,
  $C_\Davg$, and $C_\Drms$ for flow past a
  square cylinder to literature values which use finite elements. }
\end{table}

\begin{table}[H]
\center{}
  \begin{tabular}{l}
    \begin{tabular}{|c|c|c|c|c|}
      \hline
      \text{Re} &\multicolumn{2}{c|}{ L/D }  & \multicolumn{2}{c|}{$C_{\text{D}}$} \\
      \hline
      & Sen et al. & Present & Sen et al. & Present DG\\
      \hline
      5 & 2.806 & 2.786 & 4.9535 & 4.738\\
      \hline
      40 & 0.3166 & 0.323 & 1.7871 & 1.765\\
      \hline
    \end{tabular}
  \end{tabular}
  \caption{Comparison of ratio of circulation length $L$ to square side length
  $D$, and drag coefficient for steady flows.}
\end{table}
\subsection{Rear Facing Angle}
This test investigate the effects of increasingly acute angles in the interfacial geometry. We append a square cylinder of diameter $D=1$ with a rear-facing angle of varying acuteness. Six angles are considered, ranging from $\pi$ to $\frac{\pi}{9}$ (Figure 8). The square is centered at the origin of the computational domain $\Omega = [- 15, 45] \times [-30, 30]$ with $L=16$. 
For $\mathbf{u} = (u, v)$, we use incoming flow velocity $\mathbf{} u = 1$, and $v = \cos \left( \pi \frac{y}{L} \right) e^{- 2 t}$
to ensure that vortex shedding occurs at a consistent time. We set $\rho = 1$ and $\mu = \frac{1}{{Re}}$, in which $Re$ is the Reynolds number. We use
${Re} = \frac{\rho U D}{\mu}=200$ to ensure that we observe vortex
shedding. 
The effective number of grid cells on the $\ell^{\text{th}}$ level is
${n_{\ell}}  = 2^{\ell - 1} n_\cc$,  with $\ell_{\text{max}}=6$.  The time step size $\Delta t$ scales with $h_{\text{finest}}$ such that $\Delta t = \frac{
h_{\text{finest}}}{20}$, in which $h_{\text{finest}} = \frac{L}{{n_{\ell {\text{max}}}}}$. For this time step size and grid spacing, the CFL number is approximately 0.18 once the model reaches periodic steady state. We set $\kappa = 156$ and $\eta = 16$, which are determined using a bisection method. The outflow
boundary uses zero normal and tangential traction, and the top and bottom boundaries use zero tangential traction and $v = 0$.

\begin{figure}[H]
  \resizebox{1\columnwidth}{!}{\includegraphics{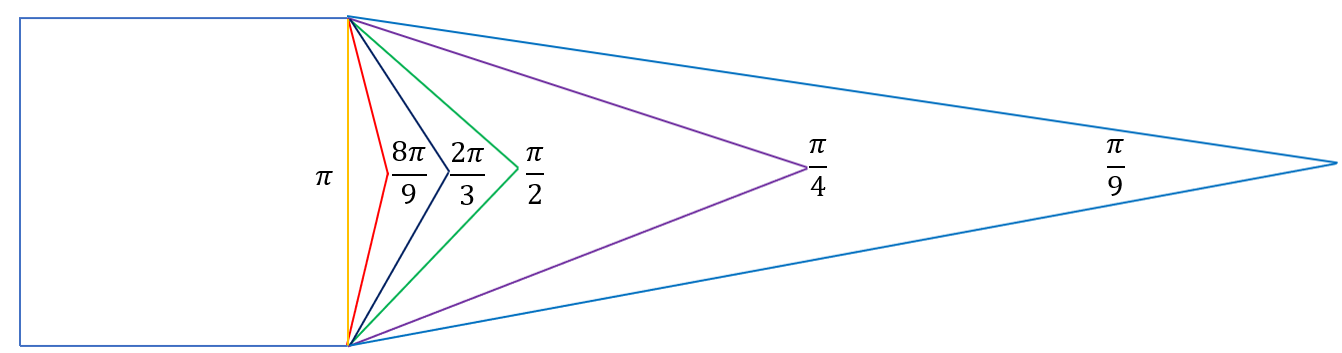}}
  
  \
  \caption{Square cylinder with rear-facing triangle, for $\theta = \left\{
  \frac{\pi}{9}, \frac{\pi}{4}, \frac{\pi}{2}, \frac{2 \pi}{3}, \frac{8
  \pi}{9}, \pi \right\}$}
\end{figure}

We vary $\theta$ in both
CG-IIM and DG-IIM formulations. Twelve cases are run with a fixed $\kappa=156$ and $\eta 
=16$ to ensure the desired accuracy is achieved. In each case, we use a bisection method
to find the largest stable $\Delta t$. For angles between $\frac{\pi}{9}$ and $\frac{8\pi}{9}$, discontinuous elements
allow a 33\% larger $\Delta t_{\max}$. The most noticable difference
occurs at $\theta = \frac{\pi}{9}$. The angle is small enough that we are unable to determine a stable $\Delta t$ for CG-IIM. In contrast, for DG-IIM, our analysis identified the same maximum stable time step size $\Delta t_\text{max}$, independent of $\theta$ (Figure 9).

\begin{figure}[H]
\center{}
  \resizebox{300pt}{!}{\includegraphics{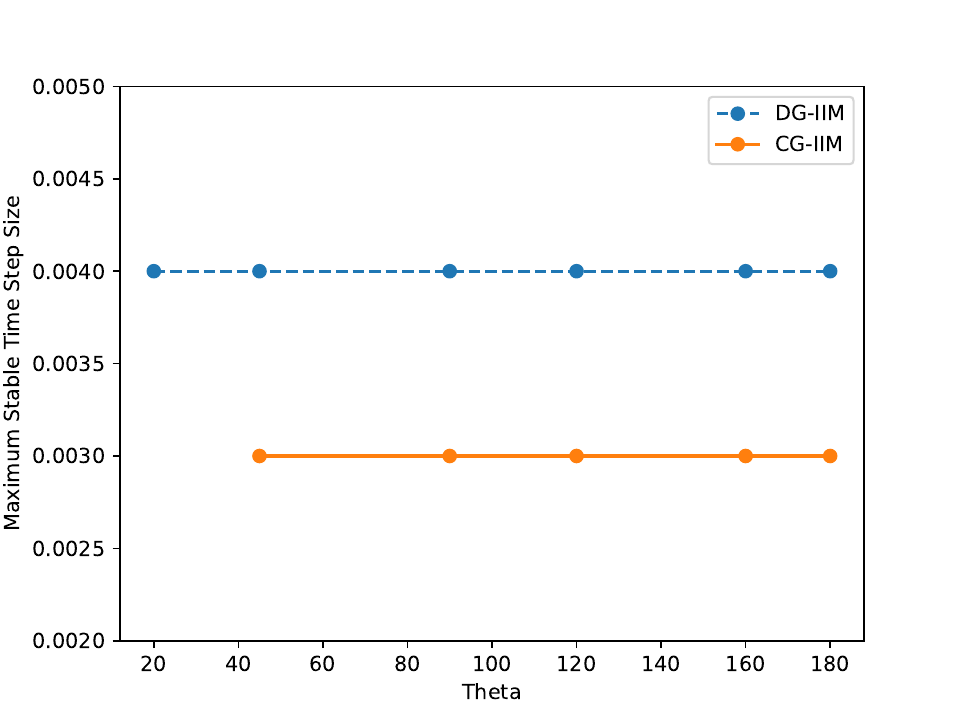   }}
  \caption{Maximum stable time step size for a square cylinder with a rear-facing
  triangle for $\kappa = 156$ and $\eta = 16$. The CG formulation for $\theta = 20$ degrees ($\frac{\pi}{9}$
  radians) could not yield a stable time step size
  with the given parameters.}
\end{figure}
As long as the angle is sufficiently large, there is no variation in stability between different degrees of acuteness for CG-IIM. For all angles, the DG formulation provides a larger $\Delta t$ for fixed $\kappa$.
\begin{figure}[H]
\center{}
  \resizebox{300pt}{!}{\includegraphics{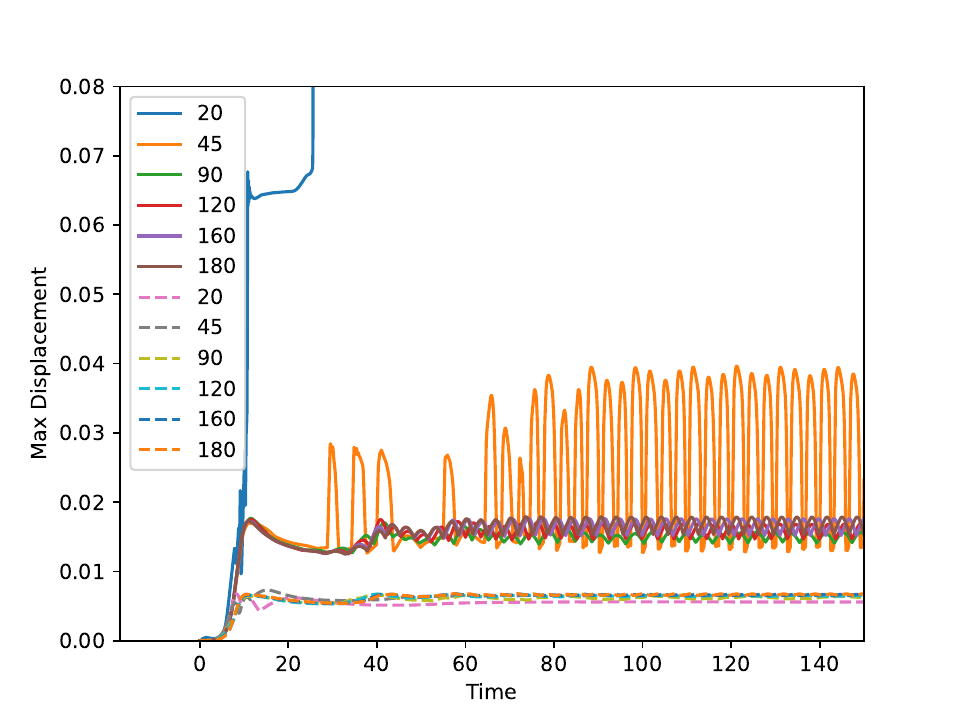}}
  
  \
  \caption{Time series of $\max (\varepsilon_{\mathbf{X}})$ for various angles (degrees). CG uses solid lines, and DG uses dashed lines. The CG formulation is not suitable for the most acute geometries.}
\end{figure}
Figure 10 shows the maximum discrepancy $\max (\varepsilon_{\mathbf{X}})$ obtained by fixing $\kappa = 156$ and $\eta = 16$. The computed values 
of $\max (\varepsilon_{\mathbf{X}})$ shows minimal variation across the angles using DG jump conditions. However, the two most acute geometries exhibit large spikes in displacement for the CG formulation. To compare the effects of the more acute geometry on the parameter sets, we vary $\kappa$ with CG jump conditions. It is possible to obtain approximately the same $\max (\varepsilon_{\mathbf{X}})$ as DG by increasing $\kappa$ by a factor of 4 (Figure 11). However, such an increase requires a proportional reduction in $\Delta t$.

\begin{figure}[H]
\center{}
  \resizebox{300pt}{!}{\includegraphics{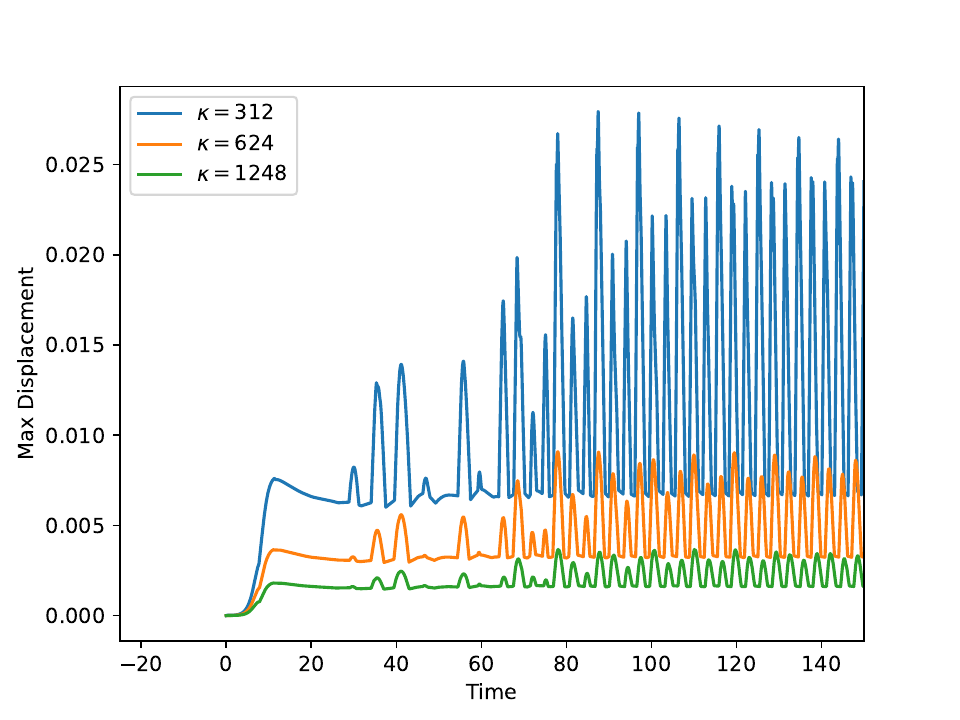}}
  
  \
  \caption{Time series of $\max (\varepsilon_{\mathbf{X}})$ for a 45 degree angle using the CG formulation. The upper bound on $\max (\varepsilon_{\mathbf{X}})$ may be reduced by increasing $\kappa$. }
\end{figure}

\subsection{Flow Past a Star}
This section considers the flow past a star-shaped interface at various Reynolds numbers. This test serves to compare the accuracy of DG-IIM and CG-IIM for a non-convex geometry. 
 The star is centered at the origin of the computational domain $\Omega = [- 8, 8]^2 $  with $L=16$. 
For $\mathbf{u} = (u, v)$, we use incoming flow velocity $\mathbf{} u = 1$, and $v = \cos \left( \pi \frac{y}{L} \right) e^{- 2 t}$
to ensure that vortex shedding occurs at a consistent time. We set $\rho = 1$ and $\mu = \frac{1}{\text{Re}}$. We use
$\text{Re} = \frac{\rho U D}{\mu}=100$, $200$, and $400$. 
The effective number of grid cells on the $\ell^{\text{th}}$ level is
${n_{\ell}}  = 2^{\ell - 1} n_\cc$,  with $\ell_{\text{max}}=6$.  The time step size $\Delta t$ is fixed at $5\cdot10^{-5}$. For this time step size and grid spacing, the CFL number is approximately 0.18 once the model reaches periodic steady state. We set $\kappa = 5\cdot10^3$ and $\eta = 16$. 
The outflow
boundary uses zero normal and tangential traction, and the top and bottom boundaries use zero tangential traction and $v = 0$.
Figure 12 shows a snapshot of the velocity field at $t=72$ after the onset of vortex shedding.

\begin{figure}[H]
\center{}
  \resizebox{200pt}{!}{\includegraphics{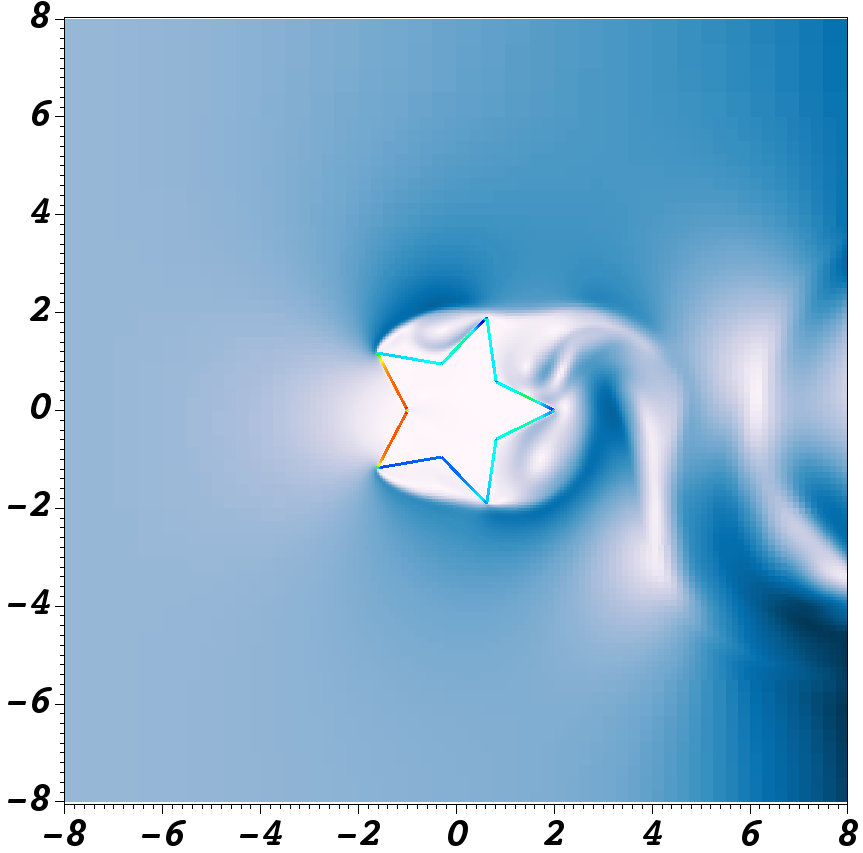}}\resizebox{80pt}{!}{

\tikzset{every picture/.style={line width=0.75pt}} 

\begin{tikzpicture}[x=0.75pt,y=0.75pt,yscale=-1,xscale=1]

\draw (83.97,256.76) node  {\includegraphics[width=50.96pt,height=79.62pt]{Figures/colorbar_displacement.png}};
\draw [color={rgb, 255:red, 0; green, 0; blue, 0 }  ,draw opacity=1 ]   (117.94,205.2) -- (141.64,205.2) ;
\draw    (117.94,309.84) -- (141.64,309.84) ;
\draw    (117.94,256.76) -- (141.64,256.76) ;
\draw [color={rgb, 255:red, 0; green, 0; blue, 0 }  ,draw opacity=1 ]   (117.94,230.98) -- (141.64,230.98) ;
\draw    (117.94,282.54) -- (141.64,282.54) ;
\draw [color={rgb, 255:red, 0; green, 0; blue, 0 }  ,draw opacity=1 ]   (119.94,60.2) -- (143.64,60.2) ;
\draw    (119.94,164.84) -- (143.64,164.84) ;
\draw    (119.94,111.76) -- (143.64,111.76) ;
\draw [color={rgb, 255:red, 0; green, 0; blue, 0 }  ,draw opacity=1 ]   (119.94,85.98) -- (143.64,85.98) ;
\draw    (119.94,137.54) -- (143.64,137.54) ;
\draw (82.94,111.58) node  {\includegraphics[width=52.5pt,height=80.07pt]{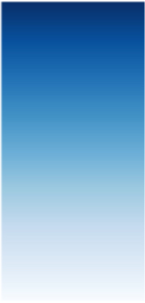}};

\draw (156.82,155.58) node [anchor=north west][inner sep=0.75pt]   [align=left] {0.00};
\draw (156.82,129.8) node [anchor=north west][inner sep=0.75pt]   [align=left] {0.62};
\draw (156.82,104.02) node [anchor=north west][inner sep=0.75pt]   [align=left] {1.22};
\draw (156.82,76.72) node [anchor=north west][inner sep=0.75pt]   [align=left] {1.84};
\draw (156.82,49.42) node [anchor=north west][inner sep=0.75pt]   [align=left] {2.45};
\draw (70,29.4) node [anchor=north west][inner sep=0.75pt]  [font=\Large]  {$||\mathbf{u}||$};
\draw (70,180) node [anchor=north west][inner sep=0.75pt]  [font=\Large]  {$\varepsilon_{\mathbf{X}}$};
\draw (154.82,194.42) node [anchor=north west][inner sep=0.75pt]   [align=left] {3.4e-4};
\draw (154.82,221.72) node [anchor=north west][inner sep=0.75pt]   [align=left] {2.5e-4};
\draw (154.82,249.02) node [anchor=north west][inner sep=0.75pt]   [align=left] {1.7e-4};
\draw (154.82,274.8) node [anchor=north west][inner sep=0.75pt]   [align=left] {9.7e-5};
\draw (154.82,300.58) node [anchor=north west][inner sep=0.75pt]   [align=left] {1.5e-5};

\end{tikzpicture}}
  \caption{Velocity field of flow past a star in the full fluid domain at $\text{Re}= 400$ using DG-IIM at $t=72$. The largest $\varepsilon_{\mathbf{X}}$ is found on the inflow-facing segments of the mesh.}
\end{figure}

We compare $\varepsilon_\mathbf{X}$ of DG-IIM and CG-IIM at $\text{Re}=100$, $200$, and $400$ to assess the two methods' impact on accuracy (Figure 13). We observe little difference in $\varepsilon_\mathbf{X}$ between the various $Re$ using DG-IIM. For all $Re$, the CG-IIM scheme results in at least three times larger $\varepsilon_\mathbf{X}$ than DG-IIM. 

\begin{figure}[H]
\center{}
  \resizebox{300pt}{!}{\includegraphics{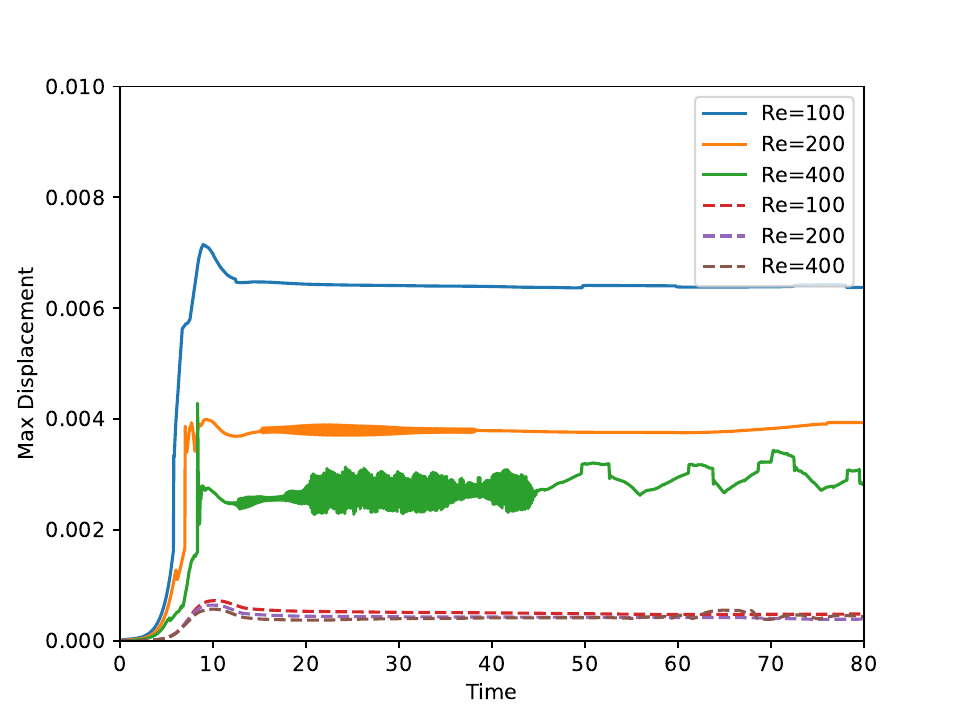}}
  \caption{Comparison of present CG-IIM and DG-IIM computations on $\varepsilon_\mathbf{X}$ over time. CG-IIM uses solid lines, while DG-IIM uses dashed lines.}
\end{figure}

\subsection{Flow Over Backward Facing Step}
This section considers the flow over a backward-facing step. This serves to demonstrate both an interior flow and a flow over non-convex geometry.

\begin{figure}[H]
\center{}
  \resizebox{300pt}{!}{

\tikzset{every picture/.style={line width=0.75pt}} 

\begin{tikzpicture}[x=0.75pt,y=0.75pt,yscale=-1,xscale=1]

\draw    (100,70) -- (388.76,70) ;
\draw    (179.76,131) -- (388.76,131) ;
\draw    (179.76,100.08) -- (179.76,131) ;
\draw    (99.76,100.08) -- (179.76,100.08) ;
\draw [color={rgb, 255:red, 74; green, 144; blue, 226 }  ,draw opacity=1 ]   (100,77) -- (117.76,77) ;
\draw [shift={(119.76,77)}, rotate = 180] [color={rgb, 255:red, 74; green, 144; blue, 226 }  ,draw opacity=1 ][line width=0.75]    (10.93,-3.29) .. controls (6.95,-1.4) and (3.31,-0.3) .. (0,0) .. controls (3.31,0.3) and (6.95,1.4) .. (10.93,3.29)   ;
\draw [color={rgb, 255:red, 74; green, 144; blue, 226 }  ,draw opacity=1 ]   (100,86) -- (124.76,86) ;
\draw [shift={(126.76,86)}, rotate = 180] [color={rgb, 255:red, 74; green, 144; blue, 226 }  ,draw opacity=1 ][line width=0.75]    (10.93,-3.29) .. controls (6.95,-1.4) and (3.31,-0.3) .. (0,0) .. controls (3.31,0.3) and (6.95,1.4) .. (10.93,3.29)   ;
\draw [color={rgb, 255:red, 74; green, 144; blue, 226 }  ,draw opacity=1 ]   (100,95) -- (117.76,95) ;
\draw [shift={(119.76,95)}, rotate = 180] [color={rgb, 255:red, 74; green, 144; blue, 226 }  ,draw opacity=1 ][line width=0.75]    (10.93,-3.29) .. controls (6.95,-1.4) and (3.31,-0.3) .. (0,0) .. controls (3.31,0.3) and (6.95,1.4) .. (10.93,3.29)   ;
\draw [color={rgb, 255:red, 74; green, 144; blue, 226 }  ,draw opacity=1 ]   (362,82) -- (379.76,82) ;
\draw [shift={(381.76,82)}, rotate = 180] [color={rgb, 255:red, 74; green, 144; blue, 226 }  ,draw opacity=1 ][line width=0.75]    (10.93,-3.29) .. controls (6.95,-1.4) and (3.31,-0.3) .. (0,0) .. controls (3.31,0.3) and (6.95,1.4) .. (10.93,3.29)   ;
\draw [color={rgb, 255:red, 74; green, 144; blue, 226 }  ,draw opacity=1 ]   (362,99) -- (386.76,99) ;
\draw [shift={(388.76,99)}, rotate = 180] [color={rgb, 255:red, 74; green, 144; blue, 226 }  ,draw opacity=1 ][line width=0.75]    (10.93,-3.29) .. controls (6.95,-1.4) and (3.31,-0.3) .. (0,0) .. controls (3.31,0.3) and (6.95,1.4) .. (10.93,3.29)   ;
\draw [color={rgb, 255:red, 74; green, 144; blue, 226 }  ,draw opacity=1 ]   (362,118) -- (379.76,118) ;
\draw [shift={(381.76,118)}, rotate = 180] [color={rgb, 255:red, 74; green, 144; blue, 226 }  ,draw opacity=1 ][line width=0.75]    (10.93,-3.29) .. controls (6.95,-1.4) and (3.31,-0.3) .. (0,0) .. controls (3.31,0.3) and (6.95,1.4) .. (10.93,3.29)   ;
\draw  [dash pattern={on 0.84pt off 2.51pt}] (100,30) -- (388.76,30) -- (388.76,192.08) -- (100,192.08) -- cycle ;

\draw (183.76,106.48) node [anchor=north west][inner sep=0.75pt]    {$h_{0}$};
\draw (273.76,135.48) node [anchor=north west][inner sep=0.75pt]    {$\ell _{o}$};
\draw (123.76,105.48) node [anchor=north west][inner sep=0.75pt]    {$\ell _{i}$};
\draw (186.76,6.48) node [anchor=north west][inner sep=0.75pt]    {$u=0,v\ =0$};
\draw (188.76,198.4) node [anchor=north west][inner sep=0.75pt]    {$u=0,v\ =0$};
\draw (54.76,41.4) node [anchor=north west][inner sep=0.75pt]    {$u=0$};
\draw (56.76,146.4) node [anchor=north west][inner sep=0.75pt]    {$u=0$};
\draw (12.76,76.4) node [anchor=north west][inner sep=0.75pt]    {$u=u_{parabolic}{}$};
\draw (394.76,42.4) node [anchor=north west][inner sep=0.75pt]    {$u=0$};
\draw (394.76,144.4) node [anchor=north west][inner sep=0.75pt]    {$u=0$};
\draw (394.76,88.4) node [anchor=north west][inner sep=0.75pt]    {$p=0$};
\end{tikzpicture}}
  \caption{Boundary condition schematic for flow over a backward-facing step.}
\end{figure}
The tube is split into two sections. In the inflow segment, the diameter is 0.52, and the diameter at the outflow segment is 1.01 such that the difference in diameters is $h_0 = 0.49$.
Inflow length $\ell_i = 5h_0 = 2.45$, and outflow length $\ell_i = 10h_0 = 9.8$. The computational domain is $\Omega = [0,12] \times [0, 3]$. The boundary conditions are shown in Figure 14. 
We set $\rho = 1$, $\mu = 0.0133$, and $\text{Re} = 100, 150, 300, 400, 450, 500, 600, 650, 700$, and $800$. We use $N= 64$ grid cells on the coarsest level and 2 levels of grid refinement with a refinement ratio of 2.
We use penalty parameters $\kappa = 8\cdot10^3$ and $\eta = 0.5$. The time step size $\Delta t$ is $4.73\cdot 10^{-5}$. 
For this time step size and grid spacing, the CFL number is approximately 0.03. Figure 15 shows a vorticity field visualization of a developed flow through the channel.
\begin{figure}[H]
\center{}
  \raisebox{0.0\height}{\includegraphics[width=10cm,height=1.5cm]{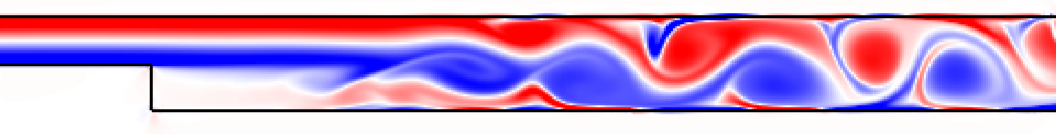}}\resizebox{50pt}{!}{

\tikzset{every picture/.style={line width=0.75pt}} 

\begin{tikzpicture}[x=0.75pt,y=0.75pt,yscale=-1,xscale=1]

\draw (78.97,110.76) node  {\includegraphics[width=50.96pt,height=79.62pt]{Figures/vorticity_colorbar.png}};
\draw [color={rgb, 255:red, 0; green, 0; blue, 0 }  ,draw opacity=1 ]   (112.94,59.2) -- (136.64,59.2) ;
\draw    (112.94,163.84) -- (136.64,163.84) ;
\draw    (112.94,110.76) -- (136.64,110.76) ;
\draw [color={rgb, 255:red, 0; green, 0; blue, 0 }  ,draw opacity=1 ]   (112.94,84.98) -- (136.64,84.98) ;
\draw    (112.94,136.54) -- (136.64,136.54) ;

\draw (70,40) node [anchor=north west][inner sep=0.75pt]  [font=\Large]  {$\omega $};
\draw (147.82,45.42) node [anchor=north west][inner sep=0.75pt]   [align=left] {150};
\draw (147.82,72.72) node [anchor=north west][inner sep=0.75pt]   [align=left] {75};
\draw (147.82,100.02) node [anchor=north west][inner sep=0.75pt]   [align=left] {0};
\draw (147.82,125.8) node [anchor=north west][inner sep=0.75pt]   [align=left] {-75};
\draw (147.82,151.58) node [anchor=north west][inner sep=0.75pt]   [align=left] {-150};

\end{tikzpicture}}

\caption{Vorticity field visualizations for developed flow over a backward facing step at $Re=800$.}

\end{figure}

We compare $\text{max}(\varepsilon_{\mathbf{X}})$ for both DG-IIM and CG-IIM over time on the lower interface. Figure 16 shows a time series comparison of the accuracy of the two methods. We observe little difference in the $\text{max}(\varepsilon_{\mathbf{X}})$  between CG-IIM and DG-IIM.

\begin{figure}[H]
\center{}
  \resizebox{300pt}{!}{\includegraphics{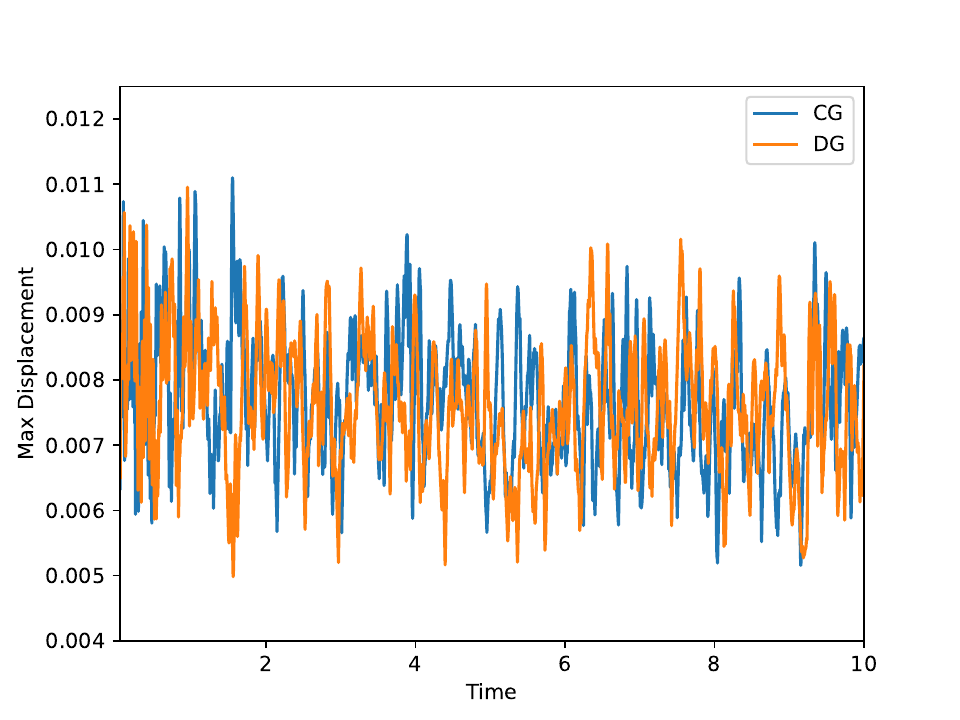}}
  \caption{Comparison of present CG-IIM and DG-IIM maximum displacements for flow over a backward facing step.}
\end{figure}

We compare the recirculation lengths across various Reynolds numbers to those reported by Armaly et al\cite{bfs} in Figure 17. Our DG-IIM results agree more closely with the authors' experimental results than the original work's two-dimensional simulations. 
\begin{figure}[H]
\center{}
  \resizebox{300pt}{!}{\includegraphics{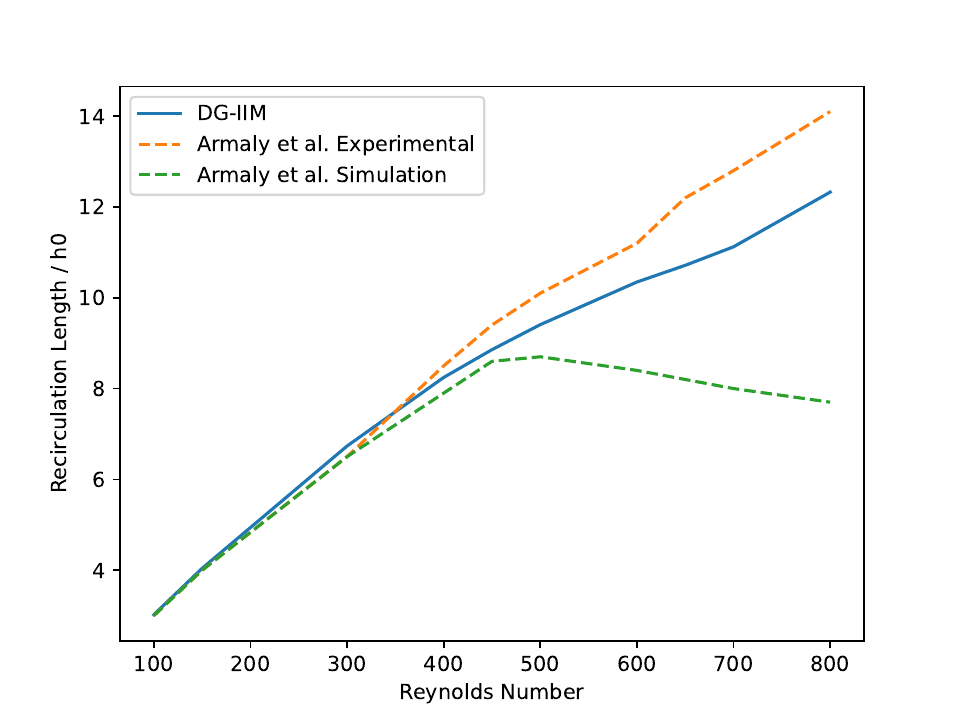}}
  \caption{Comparison of recirculation lengths for a two dimensional backward facing step at $Re = 100, 150, 300, 400, 450, 500, 600, 650, 700, $ and $800$.}
\end{figure}

\subsection{Flow Past a Sphere}

\begin{figure}[H]
\center{}
  \resizebox{300pt}{!}{\includegraphics{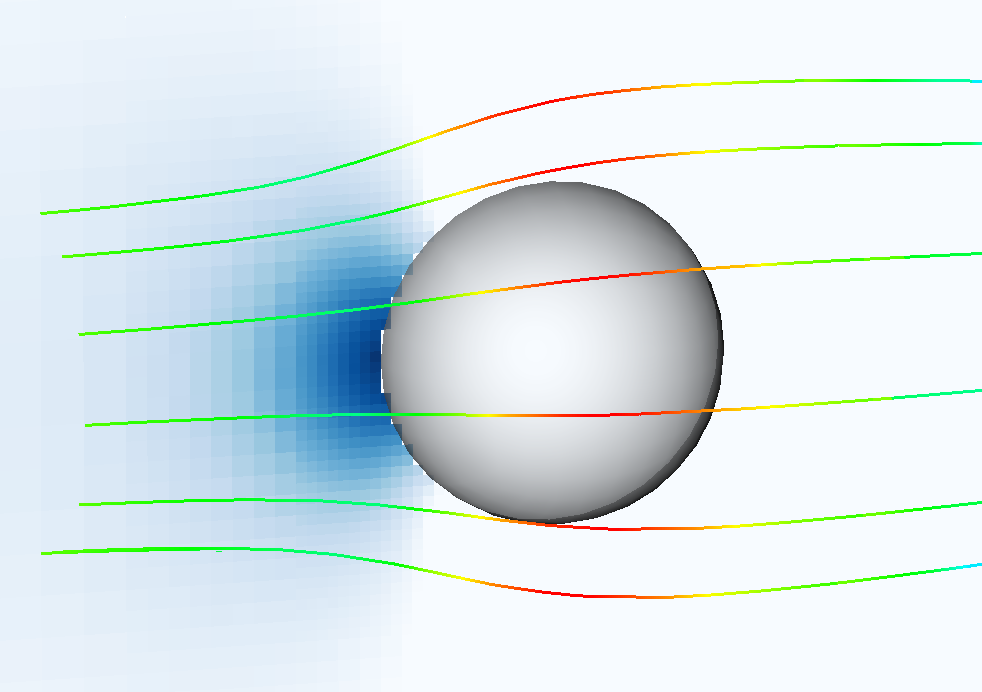}}
  \resizebox{60pt}{!}{

\tikzset{every picture/.style={line width=0.75pt}} 

\begin{tikzpicture}[x=0.75pt,y=0.75pt,yscale=-1,xscale=1]

\draw (83.97,256.76) node  {\includegraphics[width=50.96pt,height=79.62pt]{Figures/colorbar_displacement.png}};
\draw [color={rgb, 255:red, 0; green, 0; blue, 0 }  ,draw opacity=1 ]   (117.94,205.2) -- (141.64,205.2) ;
\draw    (117.94,309.84) -- (141.64,309.84) ;
\draw    (117.94,256.76) -- (141.64,256.76) ;
\draw [color={rgb, 255:red, 0; green, 0; blue, 0 }  ,draw opacity=1 ]   (117.94,230.98) -- (141.64,230.98) ;
\draw    (117.94,282.54) -- (141.64,282.54) ;
\draw [color={rgb, 255:red, 0; green, 0; blue, 0 }  ,draw opacity=1 ]   (119.94,60.2) -- (143.64,60.2) ;
\draw    (119.94,164.84) -- (143.64,164.84) ;
\draw    (119.94,111.76) -- (143.64,111.76) ;
\draw [color={rgb, 255:red, 0; green, 0; blue, 0 }  ,draw opacity=1 ]   (119.94,85.98) -- (143.64,85.98) ;
\draw    (119.94,137.54) -- (143.64,137.54) ;
\draw (82.94,111.58) node  {\includegraphics[width=52.5pt,height=80.07pt]{Figures/blue_colorbar_sphere_flow.png}};

\draw (156.82,155.58) node [anchor=north west][inner sep=0.75pt]   [align=left] {0.000};
\draw (156.82,129.8) node [anchor=north west][inner sep=0.75pt]   [align=left] {0.138};
\draw (156.82,104.02) node [anchor=north west][inner sep=0.75pt]   [align=left] {0.277};
\draw (156.82,76.72) node [anchor=north west][inner sep=0.75pt]   [align=left] {0.415};
\draw (156.82,49.42) node [anchor=north west][inner sep=0.75pt]   [align=left] {0.553};
\draw (76,29.4) node [anchor=north west][inner sep=0.75pt]  [font=\Large]  {$p$};
\draw (60,176.4) node [anchor=north west][inner sep=0.75pt]  [font=\Large]  {$||\mathbf{u} ||$};
\draw (154.82,194.42) node [anchor=north west][inner sep=0.75pt]   [align=left] {1.117};
\draw (154.82,221.72) node [anchor=north west][inner sep=0.75pt]   [align=left] {1.020};
\draw (154.82,249.02) node [anchor=north west][inner sep=0.75pt]   [align=left] {0.923};
\draw (154.82,274.8) node [anchor=north west][inner sep=0.75pt]   [align=left] {0.825};
\draw (154.82,300.58) node [anchor=north west][inner sep=0.75pt]   [align=left] {0.727};

\end{tikzpicture}}
  \caption{Pressure field slice of the $x - y$ plane of flow past a unit sphere at \text{Re}=100 using DG jump handling.
 Streamlines are colored by magnitude of velocity.}
\end{figure}

This section considers flow past a sphere (Figure 18) of diameter $D=1$ centered at the origin of the computational domain $\Omega = [- 6, 18] \times [-12, 12]\times [-12, 12] $ with $L=24$. 
For $\mathbf{u} = (u, v,w)$, we use incoming flow velocity $\mathbf{} u = 1$, $v = 0$,$w = 0$. We set $\rho = 1$ and $\mu = \frac{1}{{\text{Re}}}$, in which $Re$ is the Reynolds number. We use various Reynolds numbers to compare to data from Jones and Clarke{\cite{jones_simulation_2008}}, namely 
${\text{Re}} = \frac{\rho U D}{\mu}=40$, $60$, 
 $80$, $100$, $120$, and $140$.
The effective number of grid cells on the $\ell^{\text{th}}$ level is
${n_{\ell}}  = 2^{\ell - 1} n_\cc$,  with $\ell_{\text{max}}=7$.  The time step size $\Delta t$ scales with $h_{\text{finest}}$ such that $\Delta t = \frac{
h_{\text{finest}}}{20}$, in which $h_{\text{finest}} = \frac{L}{{n_{\ell {\text{max}}}}}$. For this time step size and grid spacing, the CFL number is approximately 0.056 once the model reaches steady state.  $\kappa = 250$ and $\eta = 1$, which are computed using a bisection method. The boundaries perpendicular to the $y$ and $z$ axes use the zero traction and no-penetration conditions. The outflow boundary uses zero normal and tangential traction.

\begin{table}[H]
\center{}
  \begin{tabular}{l}
    \begin{tabular}{|c|c|c|c|}
      \hline
      \multicolumn{4}{|c|}{$C_{\text{D}}$}  \\
      \hline
      Re & Present CG & Present DG & Jones and Clark\cite{jones_simulation_2008}\\
      \hline
      40 & 1.79 &1.77 & 1.79\\
      \hline
      60 & 1.42 & 1.41 & 1.42\\
      \hline
      80 & 1.21 &1.20 & 1.22\\
      \hline
      100 & 1.08 & 1.07 & 1.087\\
      \hline
      120 & 0.98 & 0.97 & 0.99\\
      \hline
      140 & 0.91 & 0.90 & 0.92\\
      \hline
    \end{tabular} 
  \end{tabular}
  \caption{Comparison of the steady state drag coefficient for flow past a unit sphere across various Reynolds numbers for CG-IIM and DG-IIM jump conditions and literature values.}
\end{table}

Table 5 shows agreement in $C_{\text{D}}$ between present DG-IIM computations and Jones and Clarke
{\cite{jones_simulation_2008}}. Max($\| \varepsilon_{\mathbf{X}}
\|$) was measured at ${\text{Re}} = 100$ for both CG-IIM and DG-IIM. Similar to the two
dimensional circular cylinder test, there was little difference between the two approaches. Max($\|
\varepsilon_{\mathbf{X}} \|$) was 0.0029 for DG-IIM and 0.00288 for CG-IIM.

\subsection{Flow Past a Cube}

\begin{figure}[H]
\center{}
  \resizebox{300pt}{!}{\includegraphics{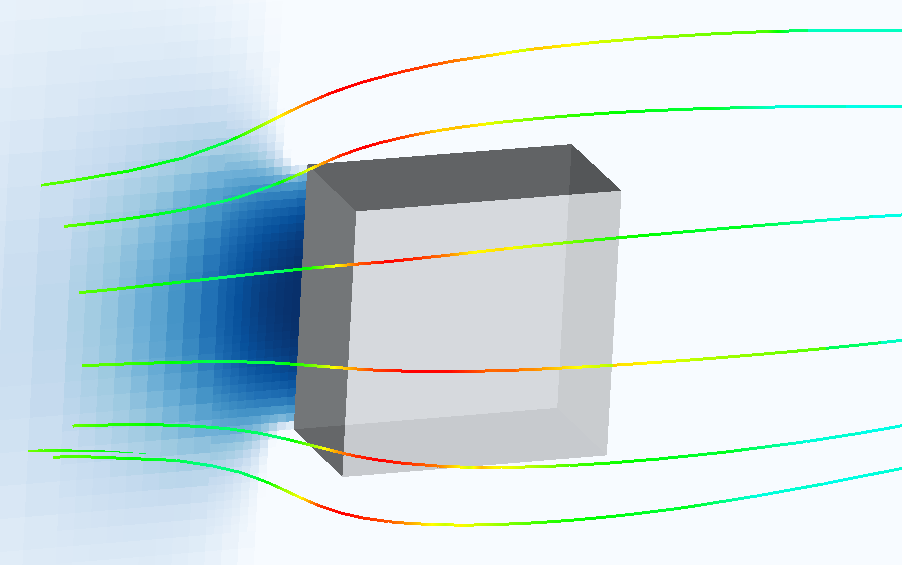}}\resizebox{60pt}{!}{

\tikzset{every picture/.style={line width=0.75pt}} 

\begin{tikzpicture}[x=0.75pt,y=0.75pt,yscale=-1,xscale=1]

\draw (83.97,256.76) node  {\includegraphics[width=50.96pt,height=79.62pt]{Figures/colorbar_displacement.png}};
\draw [color={rgb, 255:red, 0; green, 0; blue, 0 }  ,draw opacity=1 ]   (117.94,205.2) -- (141.64,205.2) ;
\draw    (117.94,309.84) -- (141.64,309.84) ;
\draw    (117.94,256.76) -- (141.64,256.76) ;
\draw [color={rgb, 255:red, 0; green, 0; blue, 0 }  ,draw opacity=1 ]   (117.94,230.98) -- (141.64,230.98) ;
\draw    (117.94,282.54) -- (141.64,282.54) ;
\draw [color={rgb, 255:red, 0; green, 0; blue, 0 }  ,draw opacity=1 ]   (119.94,60.2) -- (143.64,60.2) ;
\draw    (119.94,164.84) -- (143.64,164.84) ;
\draw    (119.94,111.76) -- (143.64,111.76) ;
\draw [color={rgb, 255:red, 0; green, 0; blue, 0 }  ,draw opacity=1 ]   (119.94,85.98) -- (143.64,85.98) ;
\draw    (119.94,137.54) -- (143.64,137.54) ;
\draw (82.94,111.58) node  {\includegraphics[width=52.5pt,height=80.07pt]{Figures/blue_colorbar_sphere_flow.png}};

\draw (156.82,155.58) node [anchor=north west][inner sep=0.75pt]   [align=left] {0.000};
\draw (156.82,129.8) node [anchor=north west][inner sep=0.75pt]   [align=left] {0.130};
\draw (156.82,104.02) node [anchor=north west][inner sep=0.75pt]   [align=left] {0.261};
\draw (156.82,76.72) node [anchor=north west][inner sep=0.75pt]   [align=left] {0.391};
\draw (156.82,49.42) node [anchor=north west][inner sep=0.75pt]   [align=left] {0.522};
\draw (76,29.4) node [anchor=north west][inner sep=0.75pt]  [font=\Large]  {$p$};
\draw (60,176.4) node [anchor=north west][inner sep=0.75pt]  [font=\Large]  {$||\mathbf{u} ||$};
\draw (154.82,194.42) node [anchor=north west][inner sep=0.75pt]   [align=left] {1.132};
\draw (154.82,221.72) node [anchor=north west][inner sep=0.75pt]   [align=left] {0.992};
\draw (154.82,249.02) node [anchor=north west][inner sep=0.75pt]   [align=left] {0.853};
\draw (154.82,274.8) node [anchor=north west][inner sep=0.75pt]   [align=left] {0.713\\};
\draw (154.82,300.58) node [anchor=north west][inner sep=0.75pt]   [align=left] {0.573};

\end{tikzpicture}}
  \caption{Pressure field slice of the $x - y$ plane of flow past a unit cube at Re=100 using DG jump handling.
 Streamlines are colored by magnitude of velocity.}
\end{figure}
This section considers flow past a cube (Figure 19) of diameter $D=1$ centered at the origin of the computational domain $\Omega = [-8.75, 26.25] \times [-17.5, 17.5]\times [-17.5, 17.5] $ with $L=35$. 
For $\mathbf{u} = (u, v,w)$, we use incoming flow velocity $\mathbf{} u = 1$, $v = 0$, $w = 0$. We set $\rho = 1$ and $\mu = \frac{1}{{Re}}$. We use various Reynolds numbers to compare to data from Saha et al.\cite{saha_three-dimensional_2004}, namely 
${\text{Re}} = \frac{\rho U D}{\mu}=50$, $100$, $150$, and $200$. 
The effective number of grid cells on the $\ell^{\text{th}}$ level is
${n_{\ell}}  = 2^{\ell - 1} n_\cc$,  with $\ell_{\text{max}}=7$.  The time step size $\Delta t$ scales with $h_{\text{finest}}$ such that $\Delta t = \frac{
h_{\text{finest}}}{10}$, in which $h_{\text{finest}} = \frac{L}{{n_{\ell {\text{max}}}}}$. For this time step size and grid spacing, the CFL number is approximately 0.03 once the model reaches steady state. We use $\kappa = 250$ and $\eta = 1$, which are determined using a bisection method. The boundaries perpendicular to the $y$ and $z$ axes use the zero traction and no-penetration conditions. The outflow boundary uses zero normal and tangential traction.

\begin{figure}[H]
\center{}
  \resizebox{300pt}{!}{\includegraphics{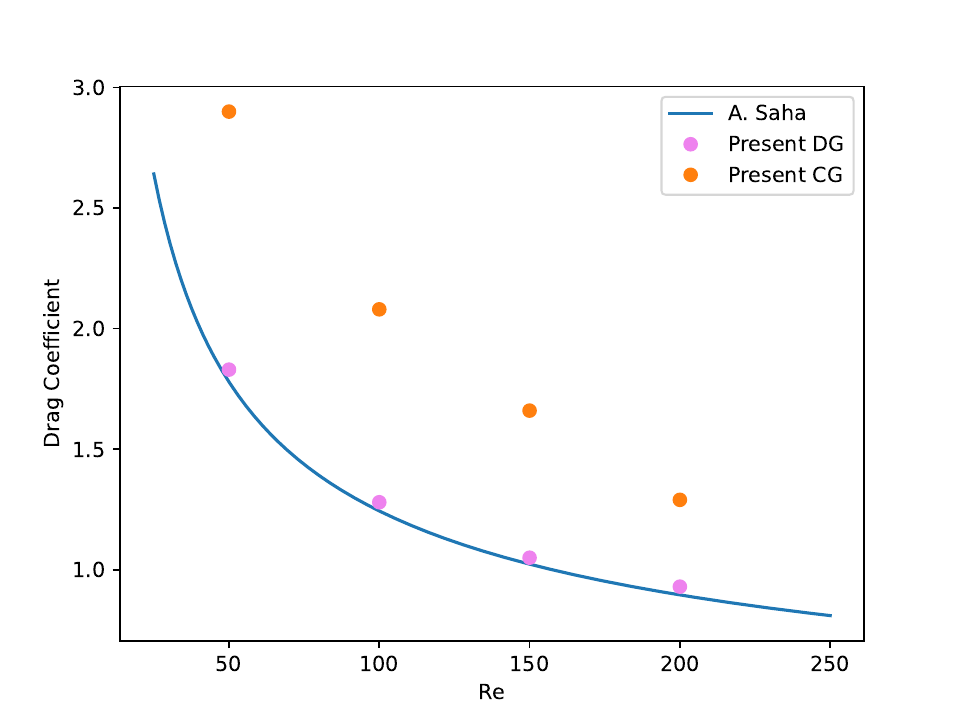}}
  \caption{Comparison of present CG-IIM and DG-IIM computations with the Reynolds-drag curve
  fit from Saha et al{\cite{saha_three-dimensional_2004}}.}
\end{figure}

Figure 20 shows the agreement between current DG-IIM simulations and data from Saha et al.\cite{saha_three-dimensional_2004} for the steady state drag coefficient at various Reynolds numbers. CG-IIM formulation overestimates  $C_{\text{D}}$ at all simulated Reynolds numbers. The
max($\varepsilon_{\mathbf{X}}$) was measured at $\text{Re} = 100$ for both
CG-IIM and DG-IIM. Similar to the two dimensional square cylinder test, there was
noticable difference for sharp geometries. The values of
max($\varepsilon_{\mathbf{X}}$) were 0.004 for DG and 0.0222 for CG.

\section{Discussion}

This study demonstrates that the IIM with DG-projected jump conditions is an effective approach to computing external flows around structures with sharp edges. The study also demonstrates that, for smooth interface geometries, the DG-IIM approach provides accuracy and efficiency that are comparable to our initial CG-IIM scheme. An important advantage of the DG-IIM representation of jump conditions is that it allows for
larger stable time steps than CG-IIM for non-smooth geometries. For the rear-facing angle model, geometry with acute angles
less than or equal to $\frac{\pi}{4}$, restrictions on $\kappa$, and as a result, $\max (\Delta t)$ become more severe for CG-IIM representations, whereas DG-IIM's maximum stable time step size need not
be decreased to maintain stability. These results suggest that DG-IIM is more
suitable to dealing with immersed interfaces with sharp features geometry, which are common in models of engineered structures. Future extensions of the
method could involve implementing higher order jump conditions, or
representing the geometry itself with a DG-IIM formulation.

\section*{Acknowledgments} 
We gratefully acknowledge research support through NSF Awards CBET 1757193, DMS 1664645, OAC 1652541, and OAC 1931516. Computations were performed using
facilities provided by University of North Carolina at Chapel Hill through the Research Computing division
of UNC Information Technology Services. We also acknowledge the help of Dr. David Wells for his mentorship in writing and debugging code, and his discussions about discontinuous basis functions. We also acknowledge Dr. Qi Sun for his helpful discussions about the IIM.


\end{document}